\documentclass[11pt,reqno]{amsart}
\usepackage[utf8]{inputenc}
\usepackage[T1]{fontenc}
\usepackage{lmodern}
\usepackage[margin=1.15in]{geometry}
\usepackage{amsmath,amssymb,mathtools,bm}
\usepackage{microtype}
\usepackage{enumitem}
\usepackage{booktabs}
\usepackage[colorlinks=true,linkcolor=blue,citecolor=blue,urlcolor=blue]{hyperref}
\numberwithin{equation}{section}
\theoremstyle{plain}
\newtheorem{theorem}{Theorem}[section]
\newtheorem{proposition}[theorem]{Proposition}
\newtheorem{lemma}[theorem]{Lemma}
\newtheorem{corollary}[theorem]{Corollary}

\newtheorem{conjecture}[theorem]{Conjecture}
\theoremstyle{definition}
\newtheorem{definition}[theorem]{Definition}
\newtheorem{assumption}[theorem]{Assumption}

\theoremstyle{remark}
\newtheorem{remark}[theorem]{Remark}
\newtheorem{example}[theorem]{Example}

\title[The circular law for non-Hermitian random band matrices]{The circular law for non-Hermitian random band matrices up to bandwidth $N^{2/3+c}$}
\author{Yi Han}
\address{Department of Mathematics, Massachusetts Institute of Technology, Cambridge, MA}
\email{hanyi16@mit.edu}
\date{}
\providecommand{\blsvP}{\mathbb P}
\providecommand{\blsvE}{\mathbb E}
\providecommand{\blsvC}{\mathbb C}
\providecommand{\blsvZ}{\mathbb Z}

\providecommand{\blsvId}{I}
\providecommand{\blsvQ}{\mathcal Q}

\providecommand{\blsvsmin}{s_{\min}}

\begin{document}
\begin{abstract}
We prove circular laws for non-Hermitian random matrices in two regimes.
First, for doubly stochastic variance profiles bounded by $C/W$ and
satisfying symmetry or transitive invariance, we establish the circular
law for real or circular complex Gaussian entries when
$W\ge N^{2/3+c}$. No lower variance bound or quantitative mixing
assumption is required. A finite-moment comparison extends the result
to bounded-density entries at explicit larger bandwidths, including
nonperiodic graph supports and inhomogeneous weights. Second, for periodic
scalar strips with a flat diagonal core, we prove a least-singular-value
estimate without a density assumption and deduce the circular law for
real subgaussian atoms, including Rademacher entries, at explicit
sublinear power bandwidths. The latter argument combines quantitative
Littlewood--Offord counting with a separator-tree elimination and does
not require full-block off-diagonal couplings. Both parts use coarse
singular-value counting to control logarithmic integrability.
\end{abstract}
\maketitle
% BEGIN sections/introduction
\section{Introduction and main results}

Let $X_N$ be an $N\times N$ matrix with eigenvalues
$\lambda_1,\ldots,\lambda_N$, counted with algebraic multiplicity.
Its empirical spectral distribution is
\[
 \mu_{X_N}:=\frac1N\sum_{j=1}^N\delta_{\lambda_j}.
\]
For a matrix with iid centered entries of variance $N^{-1}$, the
convergence of $\mu_{X_N}$ to the uniform probability measure
\[
 \mu_{\mathrm{circ}}(dz)=\frac1\pi\mathbf1_{\{|z|\le1\}}\,d^2z
\]
is the classical circular law. The Gaussian case goes back to the
ensembles introduced by Ginibre~\cite{ginibre1965}. Girko's
Hermitization method~\cite{girko1985} relates this non-Hermitian problem
to the singular values of shifted matrices. The circular law was proved
under additional regularity and moment assumptions by Bai~\cite{bai1997},
and subsequent developments include the work of
G\"otze and Tikhomirov~\cite{gotze2010} and Pan and Zhou~\cite{pan2010}.
Tao and Vu established the result under the optimal second-moment
assumption~\cite{WOS:000281425000010}. We refer to Bordenave and
Chafa\"i~\cite{bordenave2012} for a detailed account of these developments
and of the role of small singular values in the proof.

For matrices
with independent but non-identically distributed entries, the variance
profile introduces an additional dependence on the geometry of the
matrix. Random band matrices are a basic example of this setting.

In a one-dimensional band matrix, only entries within distance
$W=W_N$ of the diagonal are active. With variance normalization of order
$W^{-1}$ on the band, each row has total variance of order one.
The bandwidth interpolates between a finite-range model and a dense
random matrix. For a periodic band with sufficiently regular weights,
one expects the circular law whenever $W_N\to\infty$:
\[
 W_N\longrightarrow\infty
 \qquad\Longrightarrow\qquad
 \mu_{X_N}\Longrightarrow\mu_{\mathrm{circ}}.
\]
The question concerns the global eigenvalue distribution, and remains
meaningful when $W_N$ is much smaller than any power of $N$.
As discussed below, this threshold has been reached for several
bounded-density band models, whereas the discrete scalar-band problem
remains open at general growing bandwidth.

There is also a profile question at larger bandwidths. Write
$X_{ij}=\sqrt{s_{ij}}\,\xi_{ij}$, where the atoms are independent,
centered, and of variance one. We assume that $S=(s_{ij})$ is doubly
stochastic and that $\max_{i,j}s_{ij}\le C/W$. This normalization
includes the usual band models, but also permits holes in the support
and weights which vary with position. Even for symmetric $S$, such
profiles need not have a positive diagonal core or a translation-invariant
kernel. A circular-law theorem for a flat scalar band or an equal-weight
full-block band therefore does not directly give a theorem for this
larger class. One aim of this paper is to treat broad profiles under
density assumptions, at the cost of a larger sufficient bandwidth.

The main difficulty is to control the small singular values of
$X_N-zI_N$. Hermitization relates the empirical eigenvalue measure to
the logarithmic determinant through
\[
 \frac1N\log|\det(X_N-zI_N)|
 =\frac1N\sum_{j=1}^N\log s_j(X_N-zI_N).
\]
Weak convergence of the singular-value measure does not control the
logarithm near zero. A least-singular-value bound must be combined with
an estimate on how many singular values are small. If
$s_{\min}(X_N-zI_N)\ge e^{-L_N}$ and $0<t<1$, then
\[
 \frac1N\sum_{s_j(X_N-zI_N)\le t}|\log s_j(X_N-zI_N)|
 \le L_N\frac1N\#\{j:s_j(X_N-zI_N)\le t\}.
\]
Thus an exponential floor may suffice for the circular law, provided
that the counting estimate makes the right-hand side negligible.
Our first part obtains this control from a coarse expected resolvent
estimate and finite-moment comparison.

For discrete entries, the least-singular-value input requires a different
argument. Density-based small-ball bounds are no longer available, and
the band mask prevents a direct application of a full iid matrix
estimate. Our second part treats periodic scalar strips with a flat
diagonal core. We apply quantitative Littlewood--Offord counting to
fresh iid squares retained inside that core and propagate the resulting
bounds through the strip. This gives a density-free circular law for
real subgaussian atoms at explicit sublinear power bandwidths, including
scalar Rademacher bands.

\subsection{Previous results and scope of this paper}
The circular law has also been studied for sparse matrices whose
support is random. Wood~\cite{wood2012} proved a sparse universality
principle when the expected number of active entries in a row grows as
a fixed positive power of $N$. Basak and Rudelson~\cite{basak2019}
treated much sparser matrices with subgaussian atoms. For independent
Bernoulli selectors of parameter $p_N$ multiplying iid real atoms,
Rudelson and Tikhomirov~\cite{rudelson2019sparse} obtained the circular
law under the minimal condition $Np_N\to\infty$. In these models the
active entries are distributed throughout the matrix. A deterministic
band support imposes a different geometry, so these results do not
settle the growing-bandwidth problem for band matrices.

For independent entries with a nonconstant variance profile, the
limiting eigenvalue distribution need not be uniform on a disk.
Cook, Hachem, Najim, and Renfrew~\cite{cook2018profile,cook2022profile}
developed deterministic equivalents for such distributions and studied
their properties and examples. Their assumptions allow vanishing
variances subject to quantitative irreducibility conditions.
At the local scale, Bourgade, Yau, and Yin~\cite{MR3230002} proved a
local circular law in the iid setting, while Alt, Erd\H{o}s, and
Kr\"uger~\cite{alt2018} treated inhomogeneous profiles with uniformly
comparable variances. These results explain how the variance profile
enters the limiting law. They also indicate the additional difficulty
in allowing a sparse deterministic support without a uniform lower
bound on its nonzero variances.

Quantitative invertibility is a separate part of this problem.
Rudelson and Vershynin~\cite{rudelson2008littlewood} developed
Littlewood--Offord estimates for the smallest singular value of iid
matrices. Cook~\cite{cook2018structured} studied structured matrices
with deterministic variance profiles and additive perturbations.
Livshyts, Tikhomirov, and Vershynin~\cite{livshyts2021} obtained sharp
least-singular-value estimates for independent, non-identically
distributed entries under uniform anti-concentration assumptions.
Such assumptions cannot be imposed on the deterministic zero entries
outside a band. For discrete band matrices, our earlier
work~\cite{han2025invertibility} proves high-probability nonsingularity
for a broad class of profiles at sublinear bandwidths. The circular-law
problem additionally requires a quantitative bound for the shifted
matrix $X_N-zI_N$. In the second part of this paper, we use Jain's
combinatorial approach~\cite{blsvJain} inside the flat core of a scalar
strip to obtain such a bound.

For non-Hermitian band matrices, Jain, Jana, Luh, and
O'Rourke~\cite{jain2021circular} proved a circular law for block band
matrices with genuinely sublinear bandwidth. Tikhomirov's estimates on
the pseudospectrum of inhomogeneous matrices~\cite{tikhomirov2023pseudospectrum}
provide a least-singular-value input for more general variance
profiles. Our work~\cite{han2024circular} treats arbitrary doubly
stochastic profiles under its entry and bandwidth assumptions.
Recent work attains the optimal growing-bandwidth threshold for several
specific geometries. For open-boundary full-block tridiagonal matrices,
\cite[Theorem~1.1]{han2025slow} proves almost-sure convergence when
$N=nW$, $W\to\infty$, and $n\ge W^d$ for a sufficiently small fixed
$d>0$, under bounded-density and all-moments assumptions.
For periodic scalar bands with comparable active weights, their
polynomial endpoint tapers, and periodic full-block bands,
\cite[Theorems~2.2, 2.10 and Corollary~2.4]{han2026optimal} proves
convergence in probability for $W\to\infty$ under bounded density and
a finite third moment; a finite $(2+\alpha)$ moment also suffices with
adjusted scales. Its Theorem~2.5 treats noncompact circular-Gaussian
profiles under regularity and strict local-positivity assumptions.
Without density, its Theorem~2.8 establishes the circular law for real
subgaussian full-block matrices when $W/\log N\to\infty$.

The first part of this paper allows any symmetric doubly stochastic
variance matrix with maximum entry $O(W^{-1})$, as well as a transitive
class of nonsymmetric profiles.
For Gaussian entries, $W\ge N^{2/3+c}$ suffices, without a lower variance
bound, band-support condition, or quantitative mixing assumption.
Moment comparison extends this conclusion to bounded-density entries
at explicit larger bandwidths. This includes nonperiodic regular-graph supports and inhomogeneous
bipartite couplings which do not satisfy the profile
hypotheses of the optimal-bandwidth theorems above.

The second part removes the density assumption for scalar strips.
We obtain a least-singular-value estimate for periodic finite-range
strips having a flat diagonal core, and deduce a circular law for real
subgaussian atoms, including Rademacher entries, at explicit sublinear
power bandwidths. The off-diagonal couplings need not consist of full
iid blocks. Thus the density-free result addresses a scalar-strip
geometry outside the full-block discrete theorem just described.
For Rademacher entries, our density-free result applies when
$W\ge N^{0.93}$ for the periodic hard-cutoff indicator profile, and
when $W\ge N^{0.97}$ for the gapped-core and tapered profiles in
Examples~\ref{ex:gapped-core} and~\ref{ex:fast-taper}; see
Theorem~\ref{blsv:thm-circle} for the general exponent.
These are sufficient bandwidths, with no claim of optimality. The extensions
are in the allowed variance profiles in the first part, and in the
entry distributions for non-block scalar strips in the second.

\subsection{Correction to arXiv:2508.18143v2}
In \href{https://arxiv.org/abs/2508.18143v2}{arXiv:2508.18143v2} (26 August 2025), we claimed a circular law down to
$W\geq N^{1/2+c}$, based on a fine-scale entrywise local law for
$((X-zI)^*(X-zI)-wI)^{-1}$. The proof of that local law is incorrect:
the linearization of its self-consistent equation contains a sign error.
With the correct sign, the inverse estimate assumed in that argument does
not control the relevant operator. A nearly singular antisymmetric channel
remains, and the previous proof does not estimate the random error in its
nonconstant slow directions.

We therefore withdraw the claims that the earlier circular-law theorems
and fine local law were proved. Their statements are preserved explicitly
as conjectures and open problems in Appendix~\ref{sec:conjectures}, together
with the failed computation and the missing estimate. This is a correction
to the proof, not a counterexample to those statements. The new proved
circular-law theorem below has a different Gaussian scope and the weaker
bandwidth exponent $2/3+c$, with additional non-Gaussian extensions at
larger bandwidths. Their proofs do not use any of these conjectures or
the previous bounded-inverse assumption. Several structured subclasses of
the historical circular-law claims have meanwhile been resolved by other
arguments; these results were described above. The fine local law does not
follow from such global convergence.

\subsection{Part I: general profiles with density assumptions}
\label{subsec:general-main}
We use $[N]=\{1,\ldots,N\}$ and cyclic distance
$|i-j|_N=\min\{|i-j|,N-|i-j|\}$.

\begin{definition}[Normalized independent-entry model]\label{mainassumptionbands}
Let $S=(s_{ij})_{i,j=1}^N$ be a deterministic nonnegative matrix. Set
\[
 X_{ij}=b_{ij}\xi_{ij},\qquad b_{ij}=\sqrt{s_{ij}},
 \qquad \mathbb E\xi_{ij}=0,\quad \mathbb E|\xi_{ij}|^2=1,
\]
where the variables $\xi_{ij}$ are independent. We assume
\begin{equation}\label{eq:normalized-profile}
 \sum_j s_{ij}=\sum_j s_{ji}=1,\qquad
 s_*:=\max_{i,j}s_{ij}\leq \frac{C_0}{W}.
\end{equation}
Here $W=W_N\geq1$ and $C_0$ is independent of $N$. Thus
$W\leq C_0N$. Zero variances are allowed. 
\end{definition}

\begin{assumption}[A profile symmetry ensuring balance]\label{ass:profile-balance}
The variance profile in Definition~\ref{mainassumptionbands} satisfies at
least one of the following conditions:
\begin{enumerate}[label=\textup{(\roman*)}]
\item $S=S^T$;
\item there is a transitive subgroup $\mathfrak G_N$ of the permutations
of $[N]$ such that $s_{\pi(i),\pi(j)}=s_{ij}$ for all $\pi\in\mathfrak G_N$.
\end{enumerate}
Transitive means that any index can be sent to any other index by a member
of $\mathfrak G_N$; the group is allowed to depend on $N$.
\end{assumption}

\begin{assumption}[Gaussian entries]\label{ass:gaussian}
Assumption~\ref{ass:profile-balance} holds, and either
all $\xi_{ij}$ are standard real Gaussian variables, or all are standard
circular complex Gaussian variables $(g_{ij}+ih_{ij})/\sqrt2$, where all
$g_{ij},h_{ij}$ are independent standard real Gaussians.
\end{assumption}

Symmetry here concerns the variance profile, not the random matrix:
$X_{ij}$ and $X_{ji}$ are independent even when $s_{ij}=s_{ji}$.
The parameter $W$ measures the maximal entry variance. We do not require
that $X_{ij}$ vanish outside a band, or that positive variances satisfy a
matching lower bound.

\paragraph{Profiles motivating the first part.}
The following examples have symmetric variance profiles, but their
supports need not follow a one-dimensional band. Opposite random entries
remain independent. Here $W$ measures the inverse maximal variance,
not the geometric width of the support.

\begin{example}[Nonperiodic regular-graph supports]\label{ex:regular-graph}
Let $G_N$ be any deterministic simple undirected $d_N$-regular graph on
$[N]$, and set
\[
 s_{ij}=d_N^{-1}\mathbf1_{\{\{i,j\}\in E(G_N)\}}.
\]
Then $S=S^T$, every row and column sums to one, and $s_*=d_N^{-1}$.
No geometric ordering or expansion assumption is needed. In particular,
one may choose $G_N$ to be non-vertex-transitive; the profile then has no
transitive permutation symmetry, in any labeling. It is therefore
outside the periodic scalar and equal-weight cyclic full-block models
of \cite[Section~2]{han2026optimal}. The Gaussian theorem below applies
when $d_N\ge N^{2/3+c}$, and the subgaussian bounded-density alternative
of Theorem~\ref{thm:moment-main} applies at its stated thresholds,
with $W=d_N$.
\end{example}

\begin{example}[Inhomogeneous couplings between two populations]
\label{ex:bipartite-profile}
Let $N=2n$, $n\ge3$, and let $A$ be the $n\times n$ biadjacency matrix
of a simple $d$-regular bipartite graph. Suppose that distinct rows
$a,b$ and distinct columns $u,v$ form a four-cycle, so the four
corresponding entries of $A$ equal one. For fixed $0<\varepsilon<1/2$, set
\[
 Q=\frac1d\bigl[A+\varepsilon(e_a-e_b)(e_u-e_v)^T\bigr],
 \qquad
 S=\begin{pmatrix}0&Q\\Q^T&0\end{pmatrix}.
\]
The perturbation has zero row and column sums and preserves positivity
on the four affected edges. Thus $Q$ is doubly stochastic, $S=S^T$,
and $s_*=(1+\varepsilon)/d$. All random couplings are between the two
populations. Moreover, the squared row norms of $S$ equal
$d^{-1}+2\varepsilon^2d^{-2}$ at exactly four vertices and $d^{-1}$
at every other vertex. Hence no transitive permutation group preserves
$S$: the inhomogeneity cannot be removed by relabeling. Neither the
periodic scalar models nor the equal-weight full-block model in
\cite{han2026optimal} includes this profile. Both the Gaussian result
and the subgaussian bounded-density alternative below apply with $W=d$.
\end{example}

\begin{theorem}[Gaussian circular law]\label{thm:circular-gaussian}
Under Assumption~\ref{ass:gaussian}, fix $c\in(0,1/3)$ and suppose
$W\geq N^{2/3+c}$. Then $\mu_X$ converges weakly in probability to
$\mu_{\mathrm{circ}}$ as $N\to\infty$.
\end{theorem}

This theorem covers every symmetric doubly stochastic profile with the
indicated variance cap, as well as nonsymmetric profiles in the transitive
class. It permits reducible profiles and assumes neither expansion nor a
quantitative mixing condition. The symmetry is used only to balance two
expected diagonal resolvents. Supplementary remarks on the scope of
these assumptions appear in Section~\ref{subsec:profile-balance}.

For the non-Gaussian extension, both alternatives below impose bounded
density. The first retains the general profile scope; the second trades
the subgaussian tail assumption for a diagonal lower band.

\begin{theorem}[Moment-dependent bandwidths]\label{thm:moment-main}
Suppose Definition~\ref{mainassumptionbands} and
Assumption~\ref{ass:profile-balance} hold. Let $r\geq3$ be a fixed integer.
Assume the normalized entries match the standard real, or circular complex,
Gaussian moments through degree $r-1$, in the sense specified below, and
$\sup_{N,i,j}\mathbb E|\xi_{ij}|^r<\infty$. Suppose also that one of
the following two least-singular-value hypotheses holds:
\begin{enumerate}[label=\textup{(\roman*)}]
\item The normalized entries are uniformly subgaussian and satisfy the
bounded-density assumptions in Proposition~\ref{prop:external-lsv};
\item The bounded-density assumptions in Proposition~\ref{newleastsingularvaluebound}
hold and $s_{ij}\geq c_*/W$ whenever $|i-j|_N\leq W$, for fixed $c_*>0$.
\end{enumerate}
Then, for any fixed $c>0$ with $2(r+1)/(3r)+c<1$,
\[
 W\geq N^{\frac{2(r+1)}{3r}+c}
 \quad\Longrightarrow\quad
 \mu_X\Longrightarrow\mu_{\mathrm{circ}}\quad\text{in probability}.
\]
In case \textup{(ii)}, no moment or tail assumption beyond the displayed
finite $r$th moment is required.
\end{theorem}

For real entries, matching through degree two means centering and unit
variance; matching through degree three additionally requires
$\mathbb E\xi_{ij}^3=0$. For complex entries we match all mixed moments
of the indicated total degree to the circular Gaussian law, as specified
in \eqref{eq:moment-complex-matching}.
For real entries, the first thresholds are $8/9+c$ with a finite third
absolute moment and $5/6+c$ with zero third moment and a finite fourth
moment. The finite-$(2+\alpha)$ extension is given in
Corollary~\ref{cor:moment-finite-p}.

The main analytic input is the coarse expected count
\[
 \mathbb E\frac1N\#\{j:s_j(X-zI)\le t\}\le C_z t,
 \qquad t\ge K_0W^{-1/2},
\]
for Gaussian entries and fixed $z\ne0$. Paired with an exponential LSV
floor, this controls the logarithmic contribution of the small singular
values and gives the $2/3$ threshold. Finite-moment resolvent comparison
changes the cutoff and gives Theorem~\ref{thm:moment-main}.

\subsection{Part II: scalar strips without density assumptions}
\label{blsv:sec-statements}
In this part, $\xi$ is a fixed centered complex subgaussian random
variable with $\mathbb E|\xi|^2=1$. The LSV estimate allows complex atoms;
the circular-law conclusion is stated for real atoms.

\begin{assumption}[Periodic profile with a flat core]\label{blsv:ass-core}
Fix an integer $q\ge1$ and constants $c_0,C_0>0$. Let
$W\to\infty$, $2W+1\le N$, and let nonnegative weights
$p_{N,s}$, $-W\le s\le W$, satisfy
\begin{equation}\label{blsv:eq-profile}
 \sum_{s=-W}^W p_{N,s}=1,\qquad
 \max_s p_{N,s}\le C_0/W,\qquad
 p_{N,s}=a_N/W\quad (|s|\le\lfloor W/q\rfloor),
 \quad c_0\le a_N\le C_0.
\end{equation}
On $\blsvZ/N\blsvZ$, set
\begin{equation}\label{blsv:eq-model}
 (X_N)_{i,i+s}=\sqrt{p_{N,s}}\,\xi_{i,s},\qquad -W\le s\le W,
\end{equation}
with independent copies $\xi_{i,s}$ of $\xi$, and set the other entries
equal to zero. No symmetry $p_{N,s}=p_{N,-s}$ is required.
\end{assumption}

The variance profile is doubly stochastic and invariant under simultaneous
cyclic shifts. For $q=1$ it is exactly the flat hard-cutoff profile.
For $q=2$, only the central half-band is flat: the outer weights may vary
arbitrarily within their upper bound, including being zero. Thus this is
not a full-block random-matrix assumption on the off-diagonal couplings.

\paragraph{Profiles motivating the second part.}
Two symmetric examples show that the flat core permits genuine holes
and very small outer weights. For both, take $4W<N$ and normalize
$p_{N,s}=h_{N,s}/\sum_{|u|\le W}h_{N,u}$.

\begin{example}[A flat core with separated outer lobes]\label{ex:gapped-core}
Take
\[
 h_{N,s}=\mathbf1_{\{|s|\le\lfloor W/2\rfloor\}}
       +\mathbf1_{\{\lceil3W/4\rceil\le |s|\le W\}}.
\]
The normalization is of order $W$, and Assumption~\ref{blsv:ass-core}
holds with $q=2$. Interior holes exclude the positive-band and endpoint-taper
classes in \cite[Section~2.1]{han2026optimal}; compact support excludes
its strictly locally positive noncompact class. The support mask has
no nonzero cyclic period, so its distinct rows also rule out a full-block
profile with repeated rows within each block.
\end{example}

\begin{example}[A flat core with rapid endpoint decay]\label{ex:fast-taper}
Let $h_{N,s}=f(s/(W+1))$, where
\[
 f(x)=\begin{cases}
 1,&|x|\le1/2,\\
 \exp\!\left\{-\dfrac{(2|x|-1)^2}{(1-|x|)^2}\right\},&1/2<|x|<1,\\
 0,&|x|\ge1.
 \end{cases}
\]
Again the normalization is of order $W$ and $q=2$ is admissible.
The endpoint weight is $h_{N,W}=e^{-(W-1)^2}$, violating every fixed
polynomial endpoint lower bound in
\cite[Corollary~2.4]{han2026optimal}. The profile is compactly supported
and has unequal active weights, so neither its noncompact theorem nor
its equal-weight full-block theorem applies directly.
\end{example}

Both examples are also covered by Part I under its distributional
assumptions. The result below additionally covers Rademacher and other
real subgaussian atoms without a density. The graph-supported examples
in Part I are not asserted to satisfy this discrete theorem.

Set
\begin{equation}\label{blsv:eq-nu}
 \nu_q=q\log_2 5,\qquad
 \Theta_{N,W,q}=(N/W)^{\nu_q}\log N.
\end{equation}

\begin{theorem}[Least singular value of a discrete strip]\label{blsv:thm-strip}
Under Assumption~\ref{blsv:ass-core}, suppose
\begin{equation}\label{blsv:eq-cap-global}
 \Theta_{N,W,q}=o(W^{1/5}).
\end{equation}
For each fixed $K<\infty$, there is $C<\infty$ such that, for all
sufficiently large $N$,
\begin{equation}\label{blsv:eq-lsv-global}
 \sup_{|z|\le K}\blsvP\left\{
  \blsvsmin(X_N-z\blsvId_N)
  \le \exp\{-C\Theta_{N,W,q}\}\right\}
 \le C\frac{N}{W^{9/8}}.
\end{equation}
The same conclusion holds for $X_N+B_N-z\blsvId_N$ if $B_N$ is
deterministic, $(B_N)_{ij}=0$ for $|i-j|_N>W$, and
$\|B_N\|\le N^A$ for a fixed $A$; then $C$ may also depend on $A$.
\end{theorem}

The supremum in \eqref{blsv:eq-lsv-global} is a supremum of probabilities
for deterministic shifts. It is not a positive lower bound holding
simultaneously at all $z$ for one realization. In particular, no
arithmetic condition on $z$ is imposed. The deterministic-perturbation
assertion is an invertibility assertion only; it does not assert the
circular law for arbitrary $B_N$.

\paragraph{Comparison with Jain's deformed-iid theorem.}
Jain~\cite[Theorem~1.1]{blsvJain} considers $D+G_n$ with a full iid
noise matrix, for any centered unit-variance complex atom. At failure
probability $O(n^{-1/8})$, his stated bound gives the threshold
\[
 \bigl(C(\|D\|+\sqrt n)n^{17/8}\bigr)^{-75/2},
 \qquad \|D\|\le2^{n^{0.001}}.
\]
Under subgaussianity, Lemma~\ref{blsv:lem-local} instead gives
$(CHn^3)^{-3}$, where $H=\|D\|+C_\xi\sqrt n$, at the same failure
order, for $\|D\|\le e^{n^{1/5}}$. This reparameterizes his counting
and rounding method. Jain's theorem has the broader distributional
scope and also allows much smaller failure probabilities; the estimate
used here is tailored to repeated Schur complements. In particular,
the deformation power three leads to the recurrence factor $5^q$ and
the strip exponent $\nu_q=q\log_2 5$.
Theorem~\ref{blsv:thm-strip} then passes from full iid noise to a
scalar band mask by retaining iid squares inside each separator.
It is a geometric application of the local estimate, not a black-box
application of Jain's theorem to the full band matrix.

\begin{theorem}[Density-free circular law]\label{blsv:thm-circle}
Suppose Assumption~\ref{blsv:ass-core} holds with a fixed real centered
unit-variance subgaussian atom. Define
\[
 \alpha=\begin{cases}
  1/5,&\blsvE\xi^3=0,\\
  1/8,&\text{without an additional third-moment condition}.
 \end{cases}
 \qquad \gamma_{q,\alpha}=\frac{\nu_q}{\nu_q+\alpha}.
\]
For any fixed $\gamma\in(\gamma_{q,\alpha},1)$, if $W\ge N^\gamma$,
the empirical eigenvalue measure of $X_N$ converges weakly in probability
to the uniform probability measure on the unit disk.
\end{theorem}

For the flat hard-cutoff band ($q=1$), the circular-law threshold is
$\gamma>0.920696\ldots$ when $\mathbb E\xi^3=0$, and
$\gamma>0.948916\ldots$ for a general real subgaussian atom.
For Examples~\ref{ex:gapped-core} and~\ref{ex:fast-taper} ($q=2$),
the corresponding thresholds are $0.958711\ldots$ and
$0.973788\ldots$. In particular, $W\ge N^{0.93}$ and $W\ge N^{0.97}$
suffice for Rademacher atoms in the $q=1$ and $q=2$ cases, respectively.
These are sufficient sublinear power bandwidths; the LSV floor itself
remains exponential.

The proof combines a quantitative reparameterization of Jain's rounding
and counting argument~\cite{blsvJain} with elimination along a separator
tree. Each separator retains a fixed number of fresh iid squares in the
flat core, while the remaining couplings enter a deterministic Schur
complement. Section~\ref{blsv:sec-tree} explains the algorithm and proves the strip
and circular-law results. Section~\ref{blsv:sec-local} supplies the
technical proof of the deformed-iid estimate used there.

\paragraph{Organization of the proofs.}
Sections~\ref{sec:coarse} and~\ref{sec:circular-proof} establish the
Gaussian count and circular law. Section~\ref{sec:nongaussian} develops
moment comparison and a circular-law criterion with a separate LSV
input, then proves the bounded-density conclusions. The density-free
application is in Section~\ref{blsv:sec-tree}; its only deferred input
is proved in Section~\ref{blsv:sec-local}. Appendix~\ref{sec:conjectures}
records the former statements and the error in their proof.

% END sections/introduction
% BEGIN sections/coarse
\section{A coarse Gaussian resolvent estimate}\label{sec:coarse}

The main input is a bound on the expected regularized singular-value
density.  Unlike the local law conjectured below, this estimate does not
approximate individual resolvent entries by a deterministic solution and
does not assert overwhelming-probability control.  Its proof is carried out
at each fixed regularization scale, without a continuity argument or an
inverse stability estimate.

Throughout this section, the entries of $X=(x_{ij})_{i,j=1}^N$ are
independent centered Gaussian variables of one of the following two types:
either all entries are real Gaussian, or all entries are circular complex
Gaussian.  Zero-variance entries are allowed.  Set
\[
 s_{ij}:=\mathbb E|x_{ij}|^2,\qquad S=(s_{ij}),\qquad
 s_*:=\max_{i,j}s_{ij},
\]
and assume either that $S=S^T$, or that $S$ is invariant under a transitive
group of simultaneous index permutations, as in
Assumption~\ref{ass:profile-balance}.  No stochasticity, lower variance
bound, support condition, or spectral gap is imposed for the coarse bound.
The last subsection also records what the argument gives without either
profile condition.

For $z\in\mathbb C$ and $\eta>0$, define
\begin{equation}\label{eq:coarse-hermitization}
 Y=X-zI_N,\qquad
 H_z=\begin{pmatrix}0&Y\\Y^*&0\end{pmatrix},\qquad
 R=R_z(\eta):=(H_z-i\eta I_{2N})^{-1}.
\end{equation}
We distinguish the two copies of $[N]$ by the indices $i+$ and $i-$ and
write
\begin{equation}\label{eq:coarse-diagonal-notation}
 R_{i+,i+}=ia_i^+,\qquad R_{i-,i-}=ia_i^-,\qquad
 t_i=R_{i+,i-}.
\end{equation}
The variables $a_i^\pm$ are nonnegative, as follows from the block formulas
below.

\begin{theorem}[Coarse expected resolvent bound]\label{thm:coarse}
There are universal constants $c,C>0$ such that, under the assumptions of
this section, for every $z\ne0$ and every $\eta>0$ satisfying
\begin{equation}\label{eq:coarse-scale-condition}
 \frac{s_*}{\eta^2}\le c,
\end{equation}
one has
\begin{equation}\label{eq:coarse-coordinate-bound}
 \mathbb E a_i^+=\mathbb E a_i^-\le\frac{C}{|z|}
 \qquad (1\le i\le N).
\end{equation}
In particular,
\begin{equation}\label{eq:coarse-trace-bound}
 \mathbb E\frac1N\sum_{\alpha=1}^N
 \frac{\eta}{s_\alpha(X-zI_N)^2+\eta^2}
 \le\frac{C}{|z|}.
\end{equation}
If $s_*\le C_0/W$, these conclusions hold for every
$\eta\ge K W^{-1/2}$, where $K$ depends only on $C_0$.
More generally, at a fixed $z,\eta$ the same conclusions hold for any
Gaussian variance profile for which
$\mathbb E a_i^+=\mathbb E a_i^-$ for every $i$; the two structural
profile conditions above are sufficient ways to establish this balance.
\end{theorem}

\begin{corollary}[Coarse singular-value count]\label{cor:coarse-count}
Under the assumptions of Theorem~\ref{thm:coarse}, set
\[
 F_{N,z}(u):=\frac1N\#\{\alpha:s_\alpha(X-zI_N)\le u\}.
\]
At every scale $\eta$ satisfying \eqref{eq:coarse-scale-condition},
\begin{equation}\label{eq:coarse-count-expectation}
 \mathbb E F_{N,z}(\eta)\le C_z\eta,
 \qquad C_z:=\frac{C}{|z|},
\end{equation}
after increasing the universal constant $C$.  Consequently, for every
$u>0$,
\begin{equation}\label{eq:coarse-count-markov}
 \mathbb P\{F_{N,z}(\eta)>u\eta\}\le\frac{C_z}{u}.
\end{equation}
In particular, choosing $u=N^\varepsilon$ gives exceptional probability
at most $C_zN^{-\varepsilon}$ for each fixed $\varepsilon>0$.
\end{corollary}

\paragraph{Comparison with the counting estimate in \cite{han2024circular}.}
For an arbitrary doubly stochastic variance profile, the Gaussian
comparison of \cite[Proposition~3.2 and Corollary~3.3]{han2024circular}
controls the normalized resolvent trace, with high probability, at
heights $\eta\ge W^{-1/5}N^\varepsilon$. Combined with the LSV floor,
this gives the Gaussian bandwidth $W\ge N^{5/6+c}$ in its Theorem~1.5.
The counting argument in that paper remains a valid input for the
broader profile class.

Theorem~\ref{thm:coarse} reaches $\eta\ge KW^{-1/2}$ by a different
closure. Gaussian integration by parts and the relative variance bounds
below control the expectation directly, using
$\mathbb E a_i^+=\mathbb E a_i^-$. It suffices to have
$s_*/\eta^2$ small; no approximation by a free resolvent is required.
With the exponential LSV floor, the smaller cutoff changes the Gaussian
circular-law threshold from $5/6+c$ to $2/3+c$.
There are two distinctions in scope: we assume symmetry or transitive
invariance to obtain the pointwise expectation balance, and the result
is an expectation bound, with the probability estimate supplied by
Markov's inequality in \eqref{eq:coarse-count-markov}.
It does not replace the arbitrary-profile, high-probability comparison
of \cite{han2024circular}. For general bounded-density subgaussian
entries, the first threshold here is still $8/9+c$; the additional
moment matching in Theorem~\ref{thm:moment-main} yields the stronger
stages of the hierarchy.

\paragraph{Relation to the incorrect fine local law in arXiv:2508.18143v2.}
In \href{https://arxiv.org/abs/2508.18143v2}{arXiv:2508.18143v2}, we attempted a fine entrywise
local law and used it to claim circular-law bandwidth $N^{1/2+c}$.
Those claims are different from both the valid counting estimate of
\cite{han2024circular} and the coarse expectation bound proved here.
For clarity, the full former hypotheses and statements remain in
Appendix~\ref{sec:conjectures}, followed by the incorrect sign computation
and the missing stability estimate. The present proof does not invert
the linearized self-consistency operator whose control failed there.
It uses the relative variance estimates and a scalar positivity argument
at the coarse scale. Thus Theorem~\ref{thm:coarse} supplies a replacement
input for the revised circular-law proof; it does not validate the
former fine local law or its $1/2+c$ consequence.

\paragraph{How the estimate closes.}
We first outline the proof and introduce the quantities to be
controlled. The numbers $a_i^+$ and $a_i^-$ are diagonal entries of
positive matrices. The structural assumption on $S$ will give equal
expectations, denoted by $A_i$. We then set
\[
 B_i=\eta+\sum_j s_{ij}A_j,\qquad \tau_i=\mathbb E R_{i+,i-}.
\]
We seek an upper bound on $A_i$. The positivity of $B_i$ will allow us
to eliminate the auxiliary quantity $\tau_i$ from the expectation
equations and bound $A_i$ directly.

There are three calculations. First, Poincar\'e and Ward give relative
variance bounds of size $(s_*/\eta^2)A_i^2$. Second, Gaussian integration
by parts gives, up to correspondingly small errors,
\[
 1=A_iB_i-z\overline{\tau_i},\qquad
 0=B_i\tau_i+zA_i.
\]
Third, eliminating $\tau_i$ produces the positive expression
\[
 1\approx A_i\left(B_i+\frac{|z|^2}{B_i}\right).
\]
For any $B_i>0$, the quantity in parentheses is at least $2|z|$.
Thus a bound on $A_i$ can be obtained without estimating $B_i$ from
above and without inverting a stability operator. The calculations below
show that the errors are small relative to these positive terms when
$s_*/\eta^2$ is sufficiently small. This explains both the scale
$\eta\gtrsim W^{-1/2}$ and the form of the variance estimates we seek.

\subsection{Resolvent identities and balance of the diagonal expectations}

\begin{lemma}\label{lem:coarse-ward}
The block entries of $R$ satisfy
\begin{align}
 R_{++}&=i\eta(YY^*+\eta^2I_N)^{-1},&
 R_{--}&=i\eta(Y^*Y+\eta^2I_N)^{-1},
 \label{eq:coarse-diagonal-blocks}\\
 R_{+-}&=Y(Y^*Y+\eta^2I_N)^{-1},&
 R_{-+}&=Y^*(YY^*+\eta^2I_N)^{-1}=R_{+-}^*.
 \label{eq:coarse-offdiagonal-blocks}
\end{align}
In particular, $0\le a_i^\pm\le\eta^{-1}$ and
$R_{i-,i+}=\overline{t_i}$.  For every Hermitization index $v$,
\begin{equation}\label{eq:coarse-ward}
 \sum_u|R_{vu}|^2=\sum_u|R_{uv}|^2
 =\frac{\operatorname{Im}R_{vv}}{\eta}.
\end{equation}
Moreover,
\begin{equation}\label{eq:coarse-transpose-means}
 \mathbb E a_i^+=\mathbb E a_i^-=:A_i.
\end{equation}
For every variance profile, independently of either structural condition,
one has the pathwise identity
\begin{equation}\label{eq:coarse-exact-trace-balance}
 \sum_i a_i^+=\sum_i a_i^-.
\end{equation}
\end{lemma}

\begin{proof}
We first derive the block formulas. Since
\[
 H_z^2=\begin{pmatrix}YY^*&0\\0&Y^*Y\end{pmatrix},
 \qquad
 (H_z-i\eta I)(H_z+i\eta I)=H_z^2+\eta^2I,
\]
we have
\[
 R=(H_z+i\eta I)(H_z^2+\eta^2I)^{-1}.
\]
Multiplying these two block matrices gives
\eqref{eq:coarse-diagonal-blocks} and the expression for $R_{+-}$.
The identity
\[
 (Y^*Y+\eta^2I)^{-1}Y^*=Y^*(YY^*+\eta^2I)^{-1}
\]
follows by multiplying $(Y^*Y+\eta^2I)Y^*=Y^*(YY^*+\eta^2I)$
by the two inverses. It gives $R_{-+}=R_{+-}^*$.
The full resolvent $R$ is not Hermitian; this adjoint relation concerns
its off-diagonal blocks at the purely imaginary spectral parameter.
In particular,
\[
 a_i^+=\eta\langle e_i,(YY^*+\eta^2I)^{-1}e_i\rangle>0,
 \qquad a_i^+\le\eta\eta^{-2}=\eta^{-1},
\]
and the same argument applies to $a_i^-$.

Next, $R^*=(H_z+i\eta I)^{-1}$. The resolvent identity yields
\[
 R-R^*=R\bigl((H_z+i\eta I)-(H_z-i\eta I)\bigr)R^*
       =2i\eta RR^*.
\]
Both resolvents are functions of the same Hermitian matrix, so they
commute and $RR^*=R^*R$. Taking the $v$th diagonal entry gives
\[
 \operatorname{Im}R_{vv}
 =\eta\sum_u|R_{vu}|^2
 =\eta\sum_u|R_{uv}|^2.
\]
Thus Ward converts a sum of squared resolvent entries into a single
positive diagonal entry. We will use both the row and column versions.  The matrices $YY^*$ and $Y^*Y$ have
the same eigenvalues, which proves \eqref{eq:coarse-exact-trace-balance}.
We now explain the role of the profile condition. The equality of the
two traces does not by itself identify their diagonal entries. If
$S=S^T$, then $X^T\stackrel{d}=X$ in both Gaussian classes: transposition
permutes independent Gaussian entries with the same respective laws.
The deterministic shift is unchanged because $(zI)^T=zI$.
For deterministic $X$, the
plus block of the resolvent associated with $X^T-zI_N$ is
\[
 i\eta(Y^T\overline Y+\eta^2I_N)^{-1}
 =i\eta\bigl((Y^*Y)^T+\eta^2I_N\bigr)^{-1}.
\]
Its $i$th diagonal entry equals that of the minus block associated with
$X-zI_N$.  Taking expectations gives
\eqref{eq:coarse-transpose-means}.

Alternatively, suppose a transitive group $\mathfrak G$ leaves $S$
invariant.  For its permutation matrices $P_\pi$, independence and the
Gaussian laws give $P_\pi XP_\pi^*\stackrel{d}=X$.
Since $P_\pi(zI_N)P_\pi^*=zI_N$, the two diagonal resolvent blocks
transform by the same conjugation.  Hence
$\mathbb E a_{\pi(i)}^\pm=\mathbb E a_i^\pm$.
Transitivity makes each of these two families constant in $i$.
Taking expectations in \eqref{eq:coarse-exact-trace-balance} equates the
two constants and proves \eqref{eq:coarse-transpose-means} in this case
as well.  In the rest of the proof, only this equality is used;
no further symmetry or transitivity property of $S$ is needed.
\end{proof}

\subsection{Relative variance estimates}

For a complex random variable $Z$, we use
$\operatorname{Var}(Z)=\mathbb E|Z-\mathbb EZ|^2$.  In products of
centered variables we will explicitly specify whether complex conjugation
is present.

\begin{lemma}\label{lem:coarse-variance}
Set $\rho=4s_*/\eta^2$.  If $\rho\le1/2$, then
\begin{equation}\label{eq:coarse-relative-variance}
 \operatorname{Var}(a_i^+)+\operatorname{Var}(a_i^-)
 +\operatorname{Var}(t_i)\le C\rho A_i^2.
\end{equation}
\end{lemma}

\begin{proof}
Let $E_{uv}$ denote a matrix unit. Differentiating
$(H_z-i\eta I)R=I$ gives
\[
 (dH_z)R+(H_z-i\eta I)dR=0,\qquad dR=-R(dH_z)R.
\]
For a matrix unit, $(RE_{uv}R)_{kl}=R_{ku}R_{vl}$.
In the circular complex case, treat $x_{ab}$ and $\overline{x}_{ab}$
as separate Wirtinger coordinates. They occur at positions $(a+,b-)$
and $(b-,a+)$ of $H_z$, respectively. The derivatives are
\begin{equation}\label{eq:coarse-complex-derivatives}
 \partial_{x_{ab}}R=-RE_{a+,b-}R,\qquad
 \partial_{\overline{x}_{ab}}R=-RE_{b-,a+}R,
\end{equation}
and complex Gaussian Poincar\'e gives
\[
 \operatorname{Var}(F)\le
 \sum_{a,b}s_{ab}\mathbb E\left(
 |\partial_{x_{ab}}F|^2+|\partial_{\overline{x}_{ab}}F|^2\right).
\]
In the real Gaussian case,
\begin{equation}\label{eq:coarse-real-derivative}
 \partial_{x_{ab}}R=-R(E_{a+,b-}+E_{b-,a+})R,
\end{equation}
and real Gaussian Poincar\'e, followed by
$|u+v|^2\le2|u|^2+2|v|^2$, gives the same upper bound with an additional
factor of two.  Consequently, in either class,
\begin{align*}
 \operatorname{Var}(R_{i+,i+})
 &\le2s_*\mathbb E\sum_{a,b}
 \left(|R_{i+,a+}|^2|R_{b-,i+}|^2
       +|R_{i+,b-}|^2|R_{a+,i+}|^2\right)\\
 &\le\frac{4s_*}{\eta^2}\mathbb E(a_i^+)^2
 =\rho\mathbb E(a_i^+)^2.
\end{align*}
To see the last inequality explicitly, the first double sum factors as
\[
 \sum_{a,b}|R_{i+,a+}|^2|R_{b-,i+}|^2
 =\left(\sum_a|R_{i+,a+}|^2\right)
  \left(\sum_b|R_{b-,i+}|^2\right)
 \le\frac{(a_i^+)^2}{\eta^2}.
\]
The second double sum has the same bound. Here the sums over one half
of the indices were enlarged to sums over all Hermitization indices,
then Ward was applied. No independence of the resolvent entries is used.
Since $R_{i+,i+}=ia_i^+$, this implies
\[
 \operatorname{Var}(a_i^+)
 \le\rho\bigl(A_i^2+\operatorname{Var}(a_i^+)\bigr)
 \le\rho A_i^2+\tfrac12\operatorname{Var}(a_i^+).
\]
Equivalently,
\[
 (1-\rho)\operatorname{Var}(a_i^+)\le\rho A_i^2,
 \qquad
 \mathbb E(a_i^+)^2\le\frac{A_i^2}{1-\rho}\le2A_i^2.
\]
This yields a relative variance estimate without requiring an a priori
upper bound on $A_i$. In particular, $\operatorname{Var}(a_i^+)\le2\rho A_i^2$.
The identical argument
applies to $a_i^-$, using \eqref{eq:coarse-transpose-means}.

For the off-diagonal entry $t_i=R_{i+,i-}$, the two derivative products
are $-R_{i+,a+}R_{b-,i-}$ and $-R_{i+,b-}R_{a+,i-}$.
The sum of squared absolute values of either product over $a,b$ is at most $a_i^+a_i^-/\eta^2$, by the
row Ward identity at $i+$ and the column Ward identity at $i-$.
Thus the same derivative calculation gives
\[
 \operatorname{Var}(t_i)\le\rho\mathbb E(a_i^+a_i^-)
 \le\rho\sqrt{\mathbb E(a_i^+)^2\,\mathbb E(a_i^-)^2}
 \le2\rho A_i^2.
\]
All differentiations and Gaussian identities used here are justified at
fixed $N$ and $\eta>0$: the resolvent is smooth in the underlying real
coordinates, $\|R\|\le\eta^{-1}$, and its derivatives of any fixed order
are bounded by constants depending on that order and $\eta$.  Gaussian
variables have all moments.  Coordinates with zero variance make no
contribution to the identities.
\end{proof}

\subsection{Exact expectation equations, including the real correction}

Set
\begin{equation}\label{eq:coarse-B-tau}
 B_i:=\eta+(SA)_i=\eta+\sum_js_{ij}A_j,\qquad
 \tau_i:=\mathbb Et_i.
\end{equation}
The term $(SA)_i$ will arise by replacing the expectation of a product
of diagonal resolvent entries by the product of their expectations;
the discrepancy is a covariance, estimated using the previous lemma.
In particular $B_i\ge\eta>0$. Let $\chi=0$ in the circular complex case
and $\chi=1$ in the real case.

\begin{lemma}\label{lem:coarse-expectation-equations}
If $\rho\le1/2$, then
\begin{align}
 1&=A_iB_i-z\overline{\tau_i}+e_{1i},
 \label{eq:coarse-SD-diagonal}\\
 0&=B_i\tau_i+zA_i+e_{2i},
 \label{eq:coarse-SD-offdiagonal}
\end{align}
where
\begin{equation}\label{eq:coarse-error-bound}
 |e_{1i}|+|e_{2i}|\le C\rho A_iB_i.
\end{equation}
\end{lemma}

\begin{proof}
Write $u_{ji}=R_{j-,i+}$ and $v_{ji}=R_{j-,i-}$.  The $(i+,i+)$ and
$(i+,i-)$ entries of $(H_z-i\eta I_{2N})R=I_{2N}$ are
\begin{align}
 \sum_j(x_{ij}-z\delta_{ij})R_{j-,i+}-i\eta R_{i+,i+}&=1,
 \label{eq:coarse-raw-diagonal}\\
 \sum_j(x_{ij}-z\delta_{ij})R_{j-,i-}-i\eta R_{i+,i-}&=0.
 \label{eq:coarse-raw-offdiagonal}
\end{align}
For example, the first equation has deterministic terms
$-z\overline{t_i}+\eta a_i^+$, since
$R_{i-,i+}=\overline{t_i}$ and $-i\eta(ia_i^+)=\eta a_i^+$.
In the second equation the shift term is $-iz a_i^-$.
These factors of $i$ determine the signs after taking expectations.

Circular Gaussian integration by parts uses
$\mathbb E[x_{ij}F]=s_{ij}\mathbb E\partial_{\overline{x}_{ij}}F$.
Real Gaussian integration by parts uses
$\mathbb E[x_{ij}F]=s_{ij}\mathbb E\partial_{x_{ij}}F$ and has the extra
derivative term in \eqref{eq:coarse-real-derivative}.  The needed derivatives, with all indices displayed, are
\begin{align*}
 \partial_{\overline{x}_{ij}}R_{j-,i+}
 &=-R_{j-,j-}R_{i+,i+}=a_j^-a_i^+,\\
 \partial_{x_{ij}}R_{j-,i+}
 &=-R_{j-,i+}R_{j-,i+}=-u_{ji}^2,\\
 \partial_{\overline{x}_{ij}}R_{j-,i-}
 &=-R_{j-,j-}R_{i+,i-}=-ia_j^-t_i,\\
 \partial_{x_{ij}}R_{j-,i-}
 &=-R_{j-,i+}R_{j-,i-}=-u_{ji}v_{ji}.
\end{align*}
The notation in these four lines refers to the two formal derivative
directions. For real entries, the actual derivative is their sum.
The additional terms involve $u_{ji}^2$ and $u_{ji}v_{ji}$, which may
be complex-valued. We estimate their absolute values below using Ward's
identity. In both cases
\begin{align*}
 \mathbb E[x_{ij}R_{j-,i+}]
 &=s_{ij}\mathbb E(a_j^-a_i^+)
   -\chi s_{ij}\mathbb E u_{ji}^{2},\\
 \mathbb E[x_{ij}R_{j-,i-}]
 &=-is_{ij}\mathbb E(a_j^-t_i)
   -\chi s_{ij}\mathbb E(u_{ji}v_{ji}).
\end{align*}
Substituting these identities into the two row equations gives
\begin{align*}
 1&=\eta A_i+\sum_j s_{ij}\mathbb E(a_i^+a_j^-)
       -z\overline{\tau_i}-\chi D_{1i},\\
 0&=-i\eta\tau_i-i\sum_j s_{ij}\mathbb E(a_j^-t_i)
       -izA_i-\chi D_{2i}.
\end{align*}
Multiply the second equation by $i$, and separate the products:
\begin{align*}
 \mathbb E(a_i^+a_j^-)
 &=A_iA_j+\mathbb E[(a_i^+-A_i)(a_j^--A_j)],\\
 \mathbb E(a_j^-t_i)
 &=A_j\tau_i+\mathbb E[(a_j^--A_j)(t_i-\tau_i)].
\end{align*}
The deterministic contributions are $A_iB_i$ and $B_i\tau_i$.
This gives
\eqref{eq:coarse-SD-diagonal}--\eqref{eq:coarse-SD-offdiagonal}, with
\begin{align}
 e_{1i}
 &=\sum_js_{ij}\mathbb E[(a_i^+-A_i)(a_j^--A_j)]-\chi D_{1i},
 \label{eq:coarse-first-error}\\
 e_{2i}
 &=\sum_js_{ij}\mathbb E[(a_j^--A_j)(t_i-\tau_i)]-i\chi D_{2i},
 \label{eq:coarse-second-error}\\
 D_{1i}&:=\sum_js_{ij}\mathbb E u_{ji}^{2},\qquad
 D_{2i}:=\sum_js_{ij}\mathbb E(u_{ji}v_{ji}).
 \label{eq:coarse-real-corrections}
\end{align}
By Lemma~\ref{lem:coarse-variance} and the Cauchy--Schwarz inequality,
each sum of centered products in \eqref{eq:coarse-first-error}--\eqref{eq:coarse-second-error}
is bounded in absolute value by $C\rho A_i(SA)_i$. For example,
\[
 \left|\mathbb E[(a_j^--A_j)(t_i-\tau_i)]\right|
 \le \sqrt{\operatorname{Var}(a_j^-)\operatorname{Var}(t_i)}
 \le C\rho A_jA_i.
\]
Multiplication by $s_{ij}$ and summation in $j$ produces $(SA)_i$.
The same calculation applies to the first covariance. Keeping the
weights $s_{ij}$ in the sum gives $C\rho A_i(SA)_i$ without a
dimension-dependent factor.

For the additional terms in the real Gaussian case, Ward gives
\begin{align*}
 |D_{1i}|&\le s_*\mathbb E\sum_j|R_{j-,i+}|^2
 \le\frac{s_*}{\eta}A_i,\\
 |D_{2i}|&\le s_*\mathbb E
 \left(\sum_j|R_{j-,i+}|^2\sum_j|R_{j-,i-}|^2\right)^{1/2}
 \le\frac{s_*}{\eta}\mathbb E\sqrt{a_i^+a_i^-}
 \le\frac{s_*}{\eta}A_i.
\end{align*}
For the last inequality in the second line, Cauchy--Schwarz in
probability gives
$\mathbb E\sqrt{a_i^+a_i^-}\le
\sqrt{\mathbb Ea_i^+\,\mathbb Ea_i^-}=A_i$.
Thus the real corrections cost $s_*A_i/\eta$, whereas the covariance
terms cost $C\rho A_i(SA)_i$. Both have the required scale, since
\[
 \frac{s_*}{\eta}A_i=\frac\rho4\eta A_i\le\frac\rho4 A_iB_i,
 \qquad
 C\rho A_i(SA)_i\le C\rho A_iB_i.
\]
These estimates
prove \eqref{eq:coarse-error-bound} in both Gaussian classes.
\end{proof}

\subsection{Scalar closure and counting}

\begin{proof}[Proof of Theorem~\ref{thm:coarse}]
Since $B_i$ is real and strictly positive, the off-diagonal equation
can be solved without a stability estimate:
\[
 \tau_i=-\frac{zA_i+e_{2i}}{B_i},\qquad
 -z\overline{\tau_i}
 =\frac{|z|^2A_i}{B_i}+\frac{z\overline{e_{2i}}}{B_i}.
\]
Substitution into
\eqref{eq:coarse-SD-diagonal} yields
\begin{equation}\label{eq:coarse-scalar-equation}
 1=A_iB_i+\frac{|z|^2A_i}{B_i}+\mathcal E_i,
 \qquad
 \mathcal E_i=e_{1i}+\frac{z\overline{e_{2i}}}{B_i}.
\end{equation}
The factor $B_i$ in the second error cancels its denominator:
\[
 \left|\frac{z\overline{e_{2i}}}{B_i}\right|
 \le \frac{|z|}{B_i}\,C\rho A_iB_i=C\rho|z|A_i.
\]
Consequently \eqref{eq:coarse-error-bound} gives
\begin{equation}\label{eq:coarse-scalar-error}
 |\mathcal E_i|\le C\rho A_iB_i+C\rho|z|A_i.
\end{equation}
Since the two main terms are real and nonnegative,
\[
 1\ge A_iB_i+\frac{|z|^2A_i}{B_i}-|\mathcal E_i|
 \ge (1-C\rho)A_iB_i+\frac{|z|^2A_i}{B_i}-C\rho|z|A_i.
\]
Factoring out $A_i$ gives
\[
 A_i\left((1-C\rho)B_i+\frac{|z|^2}{B_i}-C\rho|z|\right)\le1.
\]
Choose the constant $c$ in \eqref{eq:coarse-scale-condition} sufficiently
small that $C\rho$ is a sufficiently small universal constant.  The
arithmetic--geometric mean inequality then implies
\[
 (1-C\rho)B_i+\frac{|z|^2}{B_i}-C\rho|z|
 \ge\bigl(2\sqrt{1-C\rho}-C\rho\bigr)|z|
 \ge c_1|z|
\]
for a universal $c_1>0$. This lower bound holds uniformly over
$B_i>0$, so $c_1|z|A_i\le1$. This proves
\eqref{eq:coarse-coordinate-bound}.  Finally,
\[
 \sum_i a_i^+=\eta\operatorname{Tr}(YY^*+\eta^2I_N)^{-1}
 =\sum_{\alpha=1}^N\frac{\eta}{s_\alpha(Y)^2+\eta^2},
\]
which proves \eqref{eq:coarse-trace-bound}.
\end{proof}

\begin{proof}[Proof of Corollary~\ref{cor:coarse-count}]
For $0\le s\le\eta$, one has
$\eta/(s^2+\eta^2)\ge(2\eta)^{-1}$.  Therefore
\[
 F_{N,z}(\eta)\le\frac{2\eta}{N}
 \sum_{\alpha=1}^N\frac{\eta}{s_\alpha(X-zI_N)^2+\eta^2}.
\]
Taking expectations gives $\mathbb EF_{N,z}(\eta)\le C\eta/|z|$.
For the probability statement, Markov's inequality reads
\[
 \mathbb P\{F_{N,z}(\eta)>u\eta\}
 \le\frac{\mathbb EF_{N,z}(\eta)}{u\eta}
 \le\frac{C}{u|z|}.
\]
\end{proof}

\begin{remark}\label{rem:coarse-scope}
The bound $C/|z|$ is stated for $z\ne0$. This suffices for
Hermitization, which requires convergence for almost every fixed $z$.  The proof never inverts an operator involving $S$.
In particular, slowly varying modes of a band profile do not have to be
controlled.  The loss is the coarse cutoff $\eta\asymp W^{-1/2}$ and
the expectation-only nature of \eqref{eq:coarse-trace-bound}.
The Markov estimate \eqref{eq:coarse-count-markov} should not be
interpreted as an overwhelming-probability local law with the same
threshold.
\end{remark}

% BEGIN sections/profile_balance
\subsection{Supplementary remarks on the profile assumptions}
\label{subsec:profile-balance}

The proof of Theorem~\ref{thm:coarse} and
Corollary~\ref{cor:coarse-count} is complete. This subsection discusses
the scope of the profile assumptions and is independent of the proofs
that follow. It may be skipped on a first reading; the circular-law
argument resumes in Section~\ref{sec:circular-proof}.

\paragraph{Examples covered by transitive invariance.}
Symmetry of $S$ is a sufficient condition for
\eqref{eq:coarse-transpose-means}, not a necessary condition for the
coarse estimate.  In particular, the transitive alternative genuinely
includes directed band profiles.  If
\[
 s_{ij}=q_{(j-i)\bmod N},\qquad q_k\geq0,\qquad
 \sum_{k=0}^{N-1}q_k=1,\qquad \max_kq_k\leq C_0/W,
\]
then $S$ is doubly stochastic and invariant under cyclic shifts.
No condition $q_k=q_{N-k}$ is required.  For example,
\[
 q_k=W^{-1}\mathbf1_{\{1,\ldots,W\}}(k),\qquad 2W<N,
\]
gives a one-sided band with $S\ne S^T$.
More generally, the nonzero $q_k$ need not have a matching lower bound.
Lemma~\ref{lem:coarse-ward} proves the needed balance for these profiles,
so they are covered by the coarse estimate and
Theorem~\ref{thm:circular-gaussian} with the same bandwidth exponent.

\paragraph{Support geometries that cannot be relabeled as narrow bands.}
Examples~\ref{ex:regular-graph} and~\ref{ex:bipartite-profile} illustrate that our assumptions allow more general
support geometries than a narrow band. We now explain how to ensure that
this distinction persists under every simultaneous relabeling of the
rows and columns: if the underlying graphs have uniform vertex expansion,
then every band representation has width of order $N$.

Let $G_N$ be the undirected support graph, with an edge $\{i,j\}$ whenever
$i\ne j$ and $s_{ij}>0$. Assume that, for some fixed $\kappa>0$,
\[
 |\partial U|\ge\kappa|U|\qquad(1\le |U|\le N/2),
\]
where $\partial U$ consists of the vertices outside $U$ adjacent to at
least one vertex in $U$. Label the vertices $1,\ldots,N$ in any order and arrange these labels
around a circle. The cyclic distance between labels $i$ and $j$ is
\[
 |i-j|_N=\min\{|i-j|,\,N-|i-j|\},
\]
namely the number of steps along the shorter of the two directions
around the circle. Suppose that $|i-j|_N\le h$ whenever $\{i,j\}$
is an edge.
We claim that $h\ge c_\kappa N$ for all sufficiently large $N$.
If $h\ge N/4$, this is immediate. Otherwise, take $U$ to be
$\lfloor N/2\rfloor$ consecutive vertices in this ordering. Every
neighbor outside $U$ must lie within $h$ positions of one of the two
ends of this interval. There are at most $2h$ such positions, so
\[
 \kappa\lfloor N/2\rfloor\le |\partial U|\le2h.
\]
This proves the claim and rules out bandwidth $o(N)$ in every ordering.
It also rules out a noncyclic band of width $o(N)$, since ordinary
distance bounds imply the corresponding cyclic distance bounds.

Thus, choosing graphs with this expansion property in
Examples~\ref{ex:regular-graph} and~\ref{ex:bipartite-profile} gives
profiles whose supports cannot be transformed into narrow bands by
relabeling. In Example~\ref{ex:bipartite-profile}, the four-cycle
perturbation preserves the support, so the same conclusion holds after
perturbation. Expansion is an optional condition for demonstrating this
geometric distinction; the circular-law results of Part I apply to the
examples without it.

\paragraph{A supplementary estimate without balance.}
For comparison, we record the weaker estimate obtained from the preceding
calculations when the two diagonal expectations are kept distinct.
This estimate is not used elsewhere in the paper.

\begin{proposition}[Supplementary geometric-mean estimate]
\label{prop:nonsymmetric-geometric}
Let $X$ have independent centered real Gaussian entries or independent
centered circular complex Gaussian entries, with arbitrary nonnegative
variance profile $S$.  Set
\[
 u_i=\mathbb E a_i^+,\qquad v_i=\mathbb E a_i^-,\qquad
 B_i^+=\eta+(Sv)_i,\qquad \tau_i=\mathbb E t_i.
\]
There are universal $c,C>0$ such that, if $z\ne0$ and
$s_*/\eta^2\leq c$, then
\begin{equation}\label{eq:nonsymmetric-geometric}
 \sqrt{u_iv_i}\leq\frac C{|z|}\qquad(1\leq i\leq N).
\end{equation}
Consequently, if an independent argument gives
$K^{-1}\leq u_i/v_i\leq K$ for every $i$, then
$u_i,v_i\leq C\sqrt K/|z|$, and
$\mathbb EF_{N,z}(\eta)\leq C\sqrt K\,\eta/|z|$.
\end{proposition}

\begin{proof}
All $u_i,v_i$ are positive for $\eta>0$.
Set $\rho=4s_*/\eta^2$.  The derivative and Ward computations in
Lemma~\ref{lem:coarse-variance}, without identifying the two means, give
\begin{equation}\label{eq:nonsymmetric-variances}
 \operatorname{Var}(a_i^+)\leq2\rho u_i^2,\quad
 \operatorname{Var}(a_i^-)\leq2\rho v_i^2,\quad
 \operatorname{Var}(t_i)\leq2\rho u_iv_i
\end{equation}
when $\rho\leq1/2$.
Applying the same Gaussian integration by parts to
\eqref{eq:coarse-raw-diagonal}--\eqref{eq:coarse-raw-offdiagonal} yields
\begin{equation}\label{eq:nonsymmetric-expectation-system}
 1=u_iB_i^+-z\overline{\tau_i}+e_{1i}^+,
 \qquad 0=B_i^+\tau_i+zv_i+e_{2i}^+,
\end{equation}
where
\begin{align*}
 e_{1i}^+
 &=\sum_js_{ij}\mathbb E[(a_i^+-u_i)(a_j^--v_j)]-\chi D_{1i},\\
 e_{2i}^+
 &=\sum_js_{ij}\mathbb E[(a_j^--v_j)(t_i-\tau_i)]-i\chi D_{2i}.
\end{align*}
Here $\chi,D_{1i},D_{2i}$ are as in
\eqref{eq:coarse-real-corrections}; their definitions do not use
profile symmetry.  Cauchy--Schwarz and
\eqref{eq:nonsymmetric-variances} bound the centered sums by
$C\rho u_i(Sv)_i$ and $C\rho\sqrt{u_iv_i}(Sv)_i$, respectively.
For the real corrections, Ward gives
\[
 |D_{1i}|\leq\frac{s_*}{\eta}u_i,\qquad
 |D_{2i}|\leq\frac{s_*}{\eta}\sqrt{u_iv_i}.
\]
Since $B_i^+\geq\eta$, it follows that
\begin{equation}\label{eq:nonsymmetric-errors}
 |e_{1i}^+|\leq C\rho u_iB_i^+,
 \qquad |e_{2i}^+|\leq C\rho\sqrt{u_iv_i}B_i^+.
\end{equation}
Eliminate $\tau_i$ from \eqref{eq:nonsymmetric-expectation-system}.
Taking real parts and applying \eqref{eq:nonsymmetric-errors} gives
\[
 (1-C\rho)u_iB_i^+
 +\frac{|z|^2v_i}{B_i^+}-C\rho|z|\sqrt{u_iv_i}\leq1.
\]
The first two terms have sum at least
$2\sqrt{1-C\rho}\,|z|\sqrt{u_iv_i}$ by arithmetic--geometric mean.
Choosing $\rho$ sufficiently small proves
\eqref{eq:nonsymmetric-geometric}.  The conditional conclusions follow
from $u_i=\sqrt{u_iv_i}\sqrt{u_i/v_i}$, the corresponding expression
for $v_i$, and the counting inequality used in
Corollary~\ref{cor:coarse-count}.
\end{proof}

\paragraph{Limitation of the supplementary estimate.}
The exact trace identity gives $\sum_i u_i=\sum_i v_i$ for arbitrary
$S$, but it does not convert \eqref{eq:nonsymmetric-geometric} into
an averaged resolvent bound.  Indeed,
\begin{equation}\label{eq:nonsymmetric-imbalance}
 \frac1N\sum_i u_i
 =\frac1N\sum_i\sqrt{u_iv_i}
  +\frac1{2N}\sum_i(\sqrt{u_i}-\sqrt{v_i})^2.
\end{equation}
The proposition controls the first term, not the second.
For illustration, the positive arrays
$u=(M,M^{-1})$ and $v=(M^{-1},M)$ have equal sums and
$\sqrt{u_iv_i}=1$, while their common mean tends to infinity with $M$.
These arrays are only an algebraic illustration of the insufficiency of
those two facts; they are not claimed to arise from a random-matrix model
or to disprove a nonsymmetric circular law.

Thus arbitrary doubly stochastic profiles are not covered merely by
replacing pointwise balance with trace equality.  An additional bound on
the imbalance term in \eqref{eq:nonsymmetric-imbalance}, or another
argument for the averaged resolvent, would be needed to extend this
proof to that entire class.  This leaves a gap in the available argument,
not a demonstrated necessity of profile symmetry for the result.

% END sections/profile_balance

% END sections/coarse
\section{The Gaussian circular law}\label{sec:circular-proof}
% BEGIN sections/global
\subsection{Fixed-scale singular-value convergence}\label{subsec:global}

We work at a fixed large imaginary part.  A contraction argument and
uniqueness of Stieltjes transforms supply the global singular-value input,
independently of the small-scale local law conjectured later.

\begin{proposition}\label{prop:global}
Let $X_N$ have independent centered real Gaussian entries, or independent
centered circular complex Gaussian entries, with a doubly stochastic
variance profile $S_N$ satisfying Assumption~\ref{ass:profile-balance}
(symmetric or invariant under a transitive group of simultaneous index
permutations).  Suppose
$s_{*,N}:=\max_{i,j}(S_N)_{ij}\to0$.  For every fixed $z\in\mathbb C$, the
empirical singular-value measure
\[
 \nu_{N,z}:=\frac1N\sum_{\alpha=1}^N
 \delta_{s_\alpha(X_N-zI_N)}
\]
converges weakly in probability to a deterministic probability measure
$\nu_z$ on $[0,\infty)$.  This measure is the same for both Gaussian
classes and for every such sequence of variance profiles.  In particular,
if $G_N$ is a normalized dense Ginibre matrix of either class, with
$\mathbb E|(G_N)_{ij}|^2=N^{-1}$, then for every bounded continuous
$f:[0,\infty)\to\mathbb R$,
\begin{equation}\label{eq:global-comparison}
 \int f\,d\nu_{N,z}-\int f\,d\nu_{G_N-zI_N}
 \longrightarrow0\qquad\text{in probability}.
\end{equation}
The matrices in this comparison need not be independent.
\end{proposition}

\begin{proof}
We retain the notation of Section~\ref{sec:coarse};
Lemma~\ref{lem:coarse-ward} gives the diagonal expectation balance under
either profile condition.  First fix
$\eta\in(2,3)$.  The deterministic resolvent bound and stochasticity give
\begin{equation}\label{eq:global-A-B-bounds}
 0\le A_i\le\eta^{-1},\qquad
 \eta\le B_i\le\eta+\eta^{-1}.
\end{equation}
For all sufficiently large $N$, the assumption
$\rho=4s_{*,N}/\eta^2\le1/2$ holds.  Thus the exact equations and error
bounds in Lemma~\ref{lem:coarse-expectation-equations} apply.  In
particular, \eqref{eq:coarse-scalar-equation}--\eqref{eq:coarse-scalar-error}
and \eqref{eq:global-A-B-bounds} imply
\begin{equation}\label{eq:global-approx-fixed-point}
 A_i=f_\eta((S_NA)_i)+O_z(s_{*,N}),\qquad
 f_\eta(x):=\frac{\eta+x}{(\eta+x)^2+|z|^2}.
\end{equation}
The error is uniform in $i$ and in $\eta\in(2,3)$.  Indeed,
$A_i=f_\eta((S_NA)_i)(1-\mathcal E_i)$ by
\eqref{eq:coarse-scalar-equation}, and
$|\mathcal E_i|\le C_z\,s_{*,N}$ on this fixed-scale interval.  This
argument also applies at $z=0$: it does not use the bound $C/|z|$ in
Theorem~\ref{thm:coarse}.

The map $f_\eta$ sends $[0,\eta^{-1}]$ into itself and satisfies
\[
 |f_\eta'(x)|
 =\frac{\bigl||z|^2-(\eta+x)^2\bigr|}
        {((\eta+x)^2+|z|^2)^2}
 \le\frac1{(\eta+x)^2}\le\eta^{-2}<1.
\]
Consequently there is a unique fixed point
$a_z(\eta)\in[0,\eta^{-1}]$ such that
\begin{equation}\label{eq:global-scalar-fixed-point}
 a_z(\eta)=\frac{\eta+a_z(\eta)}
 {(\eta+a_z(\eta))^2+|z|^2}.
\end{equation}
Since $S_N\mathbf1=\mathbf1$ and
$\|S_N\|_{\ell^\infty\to\ell^\infty}=1$, subtracting the constant fixed
point from \eqref{eq:global-approx-fixed-point} gives
\[
 \max_i|A_i-a_z(\eta)|
 \le\eta^{-2}\max_i|A_i-a_z(\eta)|+C_z\,s_{*,N}.
\]
It follows that
\begin{equation}\label{eq:global-mean-convergence}
 \max_i|A_i-a_z(\eta)|\le C_z\,s_{*,N}.
\end{equation}
This is a large-$\eta$ contraction only, not a stability statement down to
a vanishing imaginary part.

We also need concentration of the normalized trace.  Set
\[
 m_N(i\eta):=\frac1{2N}\operatorname{Tr}R_z(\eta)
 =i\int\frac{\eta}{s^2+\eta^2}\,d\nu_{N,z}(s).
\]
The equality follows because $YY^*$ and $Y^*Y$ have the same eigenvalues.
In the circular complex case,
\[
 \partial_{x_{ab}}m_N(i\eta)
 =-\frac{(R^2)_{b-,a+}}{2N},\qquad
 \partial_{\overline{x}_{ab}}m_N(i\eta)
 =-\frac{(R^2)_{a+,b-}}{2N}.
\]
In the real case the derivative is the sum of these two expressions.
Gaussian Poincar\'e therefore gives in both cases
\begin{equation}\label{eq:global-trace-variance}
 \operatorname{Var}(m_N(i\eta))
 \le\frac{C_{\mathrm{Poincar\acute e}}\,s_{*,N}}{N^2}\mathbb E\|R^2\|_{\mathrm{HS}}^2
 \le\frac{C_{\mathrm{Poincar\acute e}}\,s_{*,N}}{N\eta^4}.
\end{equation}
Here $C_{\mathrm{Poincar\acute e}}$ is an absolute constant, independent of $N$, $z$, $\eta$,
and the variance profile, and we recall that
$s_{*,N}=\max_{a,b}(S_N)_{ab}$. Thus $C_{\mathrm{Poincar\acute e}}\,s_{*,N}$ denotes the
product of this constant and the maximum entry variance.
The factor $s_{*,N}$ comes from bounding each variance weight in the
Gaussian Poincar\'e inequality by $s_{*,N}$.
The last inequality uses the pathwise bound
$\|R^2\|_{\mathrm{HS}}^2\le2N\eta^{-4}$.  Combining
\eqref{eq:global-mean-convergence} and
\eqref{eq:global-trace-variance} yields, for each $\eta\in(2,3)$,
\begin{equation}\label{eq:global-resolvent-test}
 \int h_\eta(s)\,d\nu_{N,z}(s)
 \longrightarrow a_z(\eta)\quad\text{in probability},
 \qquad h_\eta(s):=\frac{\eta}{s^2+\eta^2}.
\end{equation}
The same convergence holds for expectations.

We spell out why these resolvent tests determine weak convergence.
First, centering and double stochasticity imply
\begin{equation}\label{eq:global-second-moment}
 \mathbb E\int s^2\,d\nu_{N,z}(s)
 =\frac1N\mathbb E\|X_N-zI_N\|_{\mathrm{HS}}^2
 =1+|z|^2.
\end{equation}
Thus the deterministic probability measures
$\overline\nu_{N,z}:=\mathbb E\nu_{N,z}$ are tight.  Any subsequential
weak limit $\nu$ satisfies
\begin{equation}\label{eq:global-measure-characterization}
 \int h_\eta\,d\nu=a_z(\eta)
 \qquad\bigl(\eta\in\mathbb Q\cap(2,3)\bigr),
\end{equation}
because $h_\eta$ is bounded and continuous.

There is at most one probability measure on $[0,\infty)$ satisfying
\eqref{eq:global-measure-characterization}. To prove this, let $\nu_1$
and $\nu_2$ be two such measures. Define $T:[0,\infty)\to[0,\infty)$
by $T(s)=s^2$, and let $\widetilde\nu_j:=T_\#\nu_j$ be the
pushforward measure, that is,
\[
 \widetilde\nu_j(B)=\nu_j\bigl(T^{-1}(B)\bigr)
 \quad\text{for every Borel set }B\subset[0,\infty),\qquad j=1,2.
\]
Their Stieltjes transforms are
\[
 q_j(w):=\int\frac1{t-w}\,d\widetilde\nu_j(t)
       =\int\frac1{s^2-w}\,d\nu_j(s),
 \qquad w\in\mathbb C\setminus[0,\infty).
\]
For every rational $\eta\in(2,3)$,
\[
 q_j(-\eta^2)=\int\frac1{s^2+\eta^2}\,d\nu_j(s)
 =\frac1\eta\int h_\eta(s)\,d\nu_j(s)
 =\frac{a_z(\eta)}\eta,
\]
so $q_1$ and $q_2$ agree at all such points $-\eta^2$. These points
have an accumulation point inside the connected domain
$\mathbb C\setminus[0,\infty)$. The identity theorem gives $q_1=q_2$
throughout this domain, and uniqueness of the Stieltjes transform gives
$\widetilde\nu_1=\widetilde\nu_2$. Since $T$ is a homeomorphism of
$[0,\infty)$ onto itself, it follows that $\nu_1=\nu_2$.
Tightness supplies existence, so denote the resulting unique measure by
$\nu_z$.

Finally, regard $\nu_{N,z}$ as a random element of the space of probability
measures on $[0,\infty)$ with its weak topology.  For $M>0$, the set
\[
 \mathcal K_M:=\left\{\nu:\int s^2\,d\nu(s)\le M\right\}
\]
is compact: its measures are uniformly tight, and it is closed by the
lower semicontinuity of the second moment.  Markov's inequality and
\eqref{eq:global-second-moment} give
\[
 \mathbb P\{\nu_{N,z}\notin\mathcal K_M\}
 \le\frac{1+|z|^2}{M}.
\]
Hence the laws of these random measures are tight.  Every subsequential
limit in distribution satisfies
\eqref{eq:global-measure-characterization} almost surely, by
\eqref{eq:global-resolvent-test}, the continuity of integration against
each $h_\eta$, and the countability of $\mathbb Q\cap(2,3)$.  By the
uniqueness just proved, that limit equals $\nu_z$ almost surely.
Convergence in distribution to this deterministic measure is convergence
in probability.  This proves the first claim.

The normalized dense real and circular complex Ginibre matrices satisfy
the same assumptions, with $s_{*,N}=N^{-1}$.  The scalar fixed point and
the characterizing tests are identical.  Applying the result to either
Ginibre class and taking the difference of the two convergences proves
\eqref{eq:global-comparison}, regardless of the coupling of the matrices.
\end{proof}

\begin{remark}
Proposition~\ref{prop:global} proves only fixed-scale weak singular-value
convergence.  It gives neither logarithmic uniform integrability nor a
bound on the least singular value.  These are supplied separately below.
In particular, no conclusion about a non-Gaussian universality class is
being inferred from the Gaussian comparison proved here.
\end{remark}

% END sections/global
% BEGIN sections/circular
\subsection{The least-singular-value input}
We state precisely the external result needed for general Gaussian
profiles, including profiles without a diagonal lower band.

\begin{proposition}[Exponential lower bound away from $z=0$]\label{prop:external-lsv}
Suppose Definition~\ref{mainassumptionbands} holds, the normalized entries
are uniformly subgaussian, and either they are real with uniformly bounded
densities, or they are complex with independent real and imaginary parts,
each having a uniformly bounded density. Suppose $W\geq N^b$ for fixed
$b>0$. Fix $z\ne0$ and $0<\kappa<\min(3,3b/4)$. There is a deterministic
$L_N=L_N(z)\geq0$ such that
\begin{equation}\label{eq:external-lsv-event}
 \mathbb P\{s_{\min}(X-zI)<e^{-L_N}\}\leq2N^{-2},
 \qquad L_N\leq\frac{C_0}{|z|^2}\frac{N^{1+\kappa}}W
                   +|\log|z||,
\end{equation}
for all sufficiently large $N$. Neither symmetry of $S$ nor its former
bounded-inverse condition is needed here.
\end{proposition}
\begin{proof}
Apply \cite[Theorem~3.4]{tikhomirov2023pseudospectrum} with its parameter
$a=\kappa/3$, variance parameters $\sigma=1$ and $\sigma^*=\sqrt{s_*}$,
and $R\geq\max(1,|z|^{-1})$. Its lower restriction on $|z|$ is satisfied
eventually, since
\[
 \sigma^*N^{2a}\leq\sqrt{C_0}\,N^{2\kappa/3-b/2}\longrightarrow0.
\]
The resulting estimate is
\begin{equation}\label{eq:tikhomirov-exact}
 \mathbb P\left\{s_{\min}(X-zI)
 <|z|\exp\left(-\frac{N^{1+\kappa}s_*}{|z|^2}\right)\right\}
 \leq2N^{-2}.
\end{equation}
Take $L_N=C_0N^{1+\kappa}/(|z|^2W)+|\log|z||$.
The threshold $e^{-L_N}$ is no larger than the one in
\eqref{eq:tikhomirov-exact}, which proves the claim.
\end{proof}

\subsection{Logarithmic uniform integrability}
Let
\[
 \nu_{N,z}=\frac1N\sum_{j=1}^N\delta_{s_j(X-zI)}.
\]
We separate a general deterministic cutoff argument from the Gaussian
inputs.

\begin{lemma}[Two lower-tail ranges]\label{lem:log-ui}
Fix $z$ and let $\nu_{N,z}$ be random probability measures on $[0,\infty)$.
Suppose there are deterministic $\eta_0=\eta_0(N)\to0$, $L_N\geq0$,
and events $\mathcal E_N$ with $\mathbb P(\mathcal E_N^c)\to0$, such that
\begin{enumerate}[label=\textup{(\roman*)}]
\item $\mathbb E\nu_{N,z}([0,t])\leq C_z t$ for every
$\eta_0\leq t\leq1$;
\item on $\mathcal E_N$, $\nu_{N,z}([0,e^{-L_N}))=0$;
\item $\eta_0L_N\to0$.
\end{enumerate}
Then, for each $\varepsilon>0$,
\begin{equation}\label{eq:log-ui-lower}
 \lim_{\tau\downarrow0}\limsup_{N\to\infty}
 \mathbb P\left\{\int_{[0,\tau]}|\log s|\,d\nu_{N,z}(s)
                         >\varepsilon\right\}=0.
\end{equation}
Here an atom at zero makes the integral infinite.
\end{lemma}
\begin{proof}
Set $F(t)=\nu_{N,z}([0,t])$ and fix $0<\tau<1$. For $N$ large enough,
$\eta_0<\tau$. On $\mathcal E_N$,
\[
 \int_{[0,\eta_0]}|\log s|\,d\nu_{N,z}(s)\leq L_NF(\eta_0).
\]
Consequently, for every $u>0$,
\begin{equation}\label{eq:tiny-tail-probability}
 \mathbb P\left\{\int_{[0,\eta_0]}|\log s|\,d\nu_{N,z}(s)>u\right\}
 \leq\mathbb P(\mathcal E_N^c)+\frac{C_z\eta_0L_N}{u}.
\end{equation}
No independence between the counting event and $\mathcal E_N$ is used.
For the other range, Stieltjes integration by parts gives
\begin{align}
 \int_{(\eta_0,\tau]}-\log s\,d\nu_{N,z}(s)
 &=(-\log\tau)F(\tau)+(\log\eta_0)F(\eta_0)
                         +\int_{\eta_0}^{\tau}\frac{F(t)}t\,dt\notag\\
 &\leq |\log\tau|F(\tau)+\int_{\eta_0}^{\tau}\frac{F(t)}t\,dt.
 \label{eq:middle-tail-parts}
\end{align}
Taking expectations and using Tonelli's theorem,
\begin{equation}\label{eq:middle-tail-expectation}
 \mathbb E\int_{(\eta_0,\tau]}|\log s|\,d\nu_{N,z}(s)
 \leq C_z\tau(1+|\log\tau|).
\end{equation}
Apply Markov's inequality to this expression and combine it with
\eqref{eq:tiny-tail-probability}, first taking $N\to\infty$ and then
$\tau\downarrow0$.
\end{proof}

\begin{proof}[Proof of Theorem~\ref{thm:circular-gaussian}]
Fix $z\ne0$ and choose $0<\kappa<c/2$. In Proposition~\ref{prop:external-lsv}
take $b=2/3+c$; its restriction on $\kappa$ is then satisfied.
Corollary~\ref{cor:coarse-count} permits the cutoff
$\eta_0=K_0W^{-1/2}$ with fixed $K_0$ sufficiently large, and
\begin{equation}\label{eq:two-thirds-final}
 \eta_0L_N
 \leq C_z\left(\frac{N^{1+\kappa}}{W^{3/2}}+W^{-1/2}\right)
 \leq C_z\left(N^{\kappa-3c/2}+W^{-1/2}\right)\longrightarrow0.
\end{equation}
Lemma~\ref{lem:log-ui} therefore proves the lower logarithmic-tail bound.
Double stochasticity and centering give the exact second-moment identity
\begin{equation}\label{eq:second-moment-tail}
 \mathbb E\int s^2\,d\nu_{N,z}(s)
 =\frac1N\mathbb E\|X-zI\|_{\mathrm{HS}}^2=1+|z|^2.
\end{equation}
For $M\geq e$, $s\mapsto(\log s)/s^2$ decreases on $[M,\infty)$, so
\[
 \mathbb E\int_{(M,\infty)}\log s\,d\nu_{N,z}(s)
 \leq (1+|z|^2)\frac{\log M}{M^2}\longrightarrow0.
\]
This gives the upper logarithmic-tail bound, uniformly in $N$, by Markov's
inequality. In particular, no operator-norm estimate is needed here.

Let $\Gamma_N$ be a normalized dense Ginibre matrix in the same real or
complex class. All the same arguments apply to $\Gamma_N$ with $W=N$,
giving both logarithmic-tail bounds. For fixed $0<\tau<1<M$, the function
\[
 f_{\tau,M}(s)=\log\bigl(\min(M,\max(s,\tau))\bigr)
\]
is bounded and continuous on $[0,\infty)$. Proposition~\ref{prop:global}
shows that the difference of its integrals against the singular-value
measures of $X-zI$ and $\Gamma_N-zI$ tends to zero in probability.
Letting $\tau\downarrow0$ and $M\to\infty$ and using the two tail bounds,
we obtain
\begin{equation}\label{eq:replacement-log}
 \frac1N\log|\det(X-zI)|
 -\frac1N\log|\det(\Gamma_N-zI)|\longrightarrow0
 \quad\text{in probability}.
\end{equation}
Possible singular matrices are contained in events whose probabilities
tend to zero by Proposition~\ref{prop:external-lsv}.

The argument holds for every fixed $z\ne0$, hence for Lebesgue-almost
every $z\in\mathbb C$. Moreover,
$N^{-1}(\|X\|_{\mathrm{HS}}^2+\|\Gamma_N\|_{\mathrm{HS}}^2)$ is bounded
in probability by its bounded expectation. The replacement principle
\cite[Theorem~2.1]{WOS:000281425000010} now gives
$\mu_X-\mu_{\Gamma_N}\to0$ weakly in probability. The circular law for
$\Gamma_N$ completes the proof.
\end{proof}

\begin{remark}[What would improve the exponent]
The bandwidth restriction in the preceding proof arises from the
requirement $\eta_0L_N\to0$; the macroscopic comparison imposes no
additional restriction. A count of order $N^{1+\varepsilon}/W$ at singular-value scale
$W^{-1}$, as in the former local-law claim, would instead lead to the
product $N^{1+\kappa+\varepsilon}/W^2$. Such a count is not proved here.
Thus the $1/2+c$ assertion is not a consequence of the present coarse
estimate merely by changing parameters.
\end{remark}

% END sections/circular
% BEGIN sections/conjectures
\section{Non-Gaussian entries and a bandwidth hierarchy}\label{sec:nongaussian}

The Gaussian coarse estimate can also be used for non-Gaussian matrices,
at the cost of a larger regularization scale.  We give an elementary
comparison with a Gaussian matrix having the same variance profile.
The comparison requires only a finite absolute moment and matching of
the lower moments.  Least-singular-value estimates remain a separate
input: moment matching alone does not supply them.

\subsection{Moment assumptions and the comparison estimate}

Let $X$ satisfy Definition~\ref{mainassumptionbands}, and let
$X^{\mathrm G}_{ij}=b_{ij}g_{ij}$ be its Gaussian reference matrix, with
independent variables $g_{ij}$, also independent of all the $\xi_{ij}$.
In the real case, all $\xi_{ij}$ are real and the $g_{ij}$ are standard
real Gaussians.  In the complex case, all $g_{ij}$ are standard circular
complex Gaussians, and we require
\begin{equation}\label{eq:moment-complex-covariance}
 \mathbb E\xi_{ij}^2=0
\end{equation}
in addition to centering and $\mathbb E|\xi_{ij}|^2=1$.  Thus the real
two-dimensional covariance of each complex entry agrees with that of
its reference Gaussian.  Independence of the real and imaginary parts
is not needed for the comparison below.

\begin{definition}[Matching through degree $r-1$]\label{def:moment-matching}
For an integer $r\ge3$, we say the entries satisfy $\mathcal M_r$ if
there is a constant $M_r<\infty$, independent of $N,i,j$, such that
\[
 \mathbb E|\xi_{ij}|^r\le M_r,
\]
and the following matching conditions hold.  In the real case,
\begin{equation}\label{eq:moment-real-matching}
 \mathbb E\xi_{ij}^k=\mathbb Eg_{ij}^k,
 \qquad 0\le k\le r-1.
\end{equation}
In the complex case,
\begin{equation}\label{eq:moment-complex-matching}
 \mathbb E\bigl[\xi_{ij}^a\overline{\xi_{ij}}^{b}\bigr]
 =\mathbb E\bigl[g_{ij}^a\overline{g_{ij}}^{b}\bigr],
 \qquad a,b\ge0,\quad a+b\le r-1.
\end{equation}
For a standard circular complex Gaussian, the right-hand side equals
$\mathbf1_{\{a=b\}}a!$.
\end{definition}

For $\eta>0$ and fixed $z\in\mathbb C$, write
\[
 m_X(i\eta)=\frac1{2N}\operatorname{Tr}(H_z(X)-i\eta I_{2N})^{-1},
 \qquad H_z(X)=\begin{pmatrix}0&X-zI_N\\(X-zI_N)^*&0\end{pmatrix}.
\]

\begin{proposition}[Finite-moment resolvent comparison]
\label{prop:moment-comparison}
Suppose $\mathcal M_r$ holds for a fixed integer $r\ge3$.  Then, for every
$\eta>0$ and every $z\in\mathbb C$,
\begin{equation}\label{eq:moment-comparison}
 \left|\mathbb Em_X(i\eta)-\mathbb Em_{X^{\mathrm G}}(i\eta)\right|
 \le \frac{C_r(1+M_r)s_*^{(r-2)/2}}{\eta^{r+1}}.
\end{equation}
This comparison does not require the profile balance condition in
Assumption~\ref{ass:profile-balance}; that condition will be used when
applying the Gaussian coarse bound to $X^{\mathrm G}$.
\end{proposition}

\begin{proof}
Replace the independent entries one at a time.  Conditional on all
entries other than the one being replaced, let $H_0$ be the Hermitization
in which that random entry is zero; the deterministic shift $-zI_N$ is
unchanged.  If the entry is in position $(i,j)$, set
\[
 V(x)=xE_{i+,j-}+\overline xE_{j-,i+},\qquad
 R_0=(H_0-i\eta I_{2N})^{-1},\qquad
 R_x=(H_0+V(x)-i\eta I_{2N})^{-1}.
\]
Then $V(x)$ is Hermitian, has rank at most two, and satisfies
$\|V(x)\|=|x|$.  The resolvent identity, iterated a finite number of
times, gives the exact expansion
\begin{equation}\label{eq:moment-finite-expansion}
 R_x=\sum_{k=0}^{r-1}(-1)^k(R_0V(x))^kR_0
       +(-1)^r(R_0V(x))^rR_x.
\end{equation}
There is no assumption that $|x|/\eta$ is small.  All resolvents in
\eqref{eq:moment-finite-expansion} have norm at most $\eta^{-1}$.
The remainder has rank at most two, so
\begin{equation}\label{eq:moment-remainder-rank}
 \left|\frac1{2N}\operatorname{Tr}(R_0V(x))^rR_x\right|
 \le\frac{|x|^r}{N\eta^{r+1}}.
\end{equation}

The normalized trace of the $k$th term in the finite sum is a homogeneous
polynomial of degree $k$ in $x$ in the real case, and in $x,\overline x$
in the complex case.  Because $H_0$ is independent of the entry being
replaced, the conditional expectations of all these terms cancel under
\eqref{eq:moment-real-matching} or
\eqref{eq:moment-complex-matching}.  Applying
\eqref{eq:moment-remainder-rank} to $x=b_{ij}\xi_{ij}$ and
$x=b_{ij}g_{ij}$ bounds the contribution of this replacement by
\[
 \frac{s_{ij}^{r/2}}{N\eta^{r+1}}
       \bigl(\mathbb E|\xi_{ij}|^r+\mathbb E|g_{ij}|^r\bigr).
\]
Finally, double stochasticity implies $\sum_{i,j}s_{ij}=N$, and hence
\[
 \frac1N\sum_{i,j}s_{ij}^{r/2}
 \le s_*^{(r-2)/2}.
\]
Summing over all replacements proves \eqref{eq:moment-comparison}.
\end{proof}

\begin{remark}[What a vanishing third moment means]\label{rem:complex-third}
For real entries, $\mathcal M_3$ requires no matching beyond centering
and variance one.  If the fourth absolute moments are uniformly bounded,
the additional condition $\mathbb E\xi_{ij}^3=0$ gives $\mathcal M_4$.
For complex entries, $\mathbb E\xi_{ij}^3=0$ alone is insufficient.
In addition to \eqref{eq:moment-complex-covariance}, matching all third
moments requires
\begin{equation}\label{eq:complex-third-matching-explicit}
 \mathbb E\xi_{ij}^3=0,\qquad
 \mathbb E\bigl[\xi_{ij}^2\overline{\xi_{ij}}\bigr]=0,
\end{equation}
and their complex conjugates.  Central symmetry
$\xi_{ij}\stackrel{d}=-\xi_{ij}$ implies all these odd-degree identities.
No symmetry or identical-distribution condition between $\xi_{ij}$ and
$\xi_{ji}$ is otherwise needed: the profile condition is imposed on
the Gaussian reference, not on the full non-Gaussian laws.
\end{remark}

\subsection{An extension with only a finite moment above two}

\begin{proposition}\label{prop:moment-fractional}
Let $2<p\le3$.  Suppose the entries match the first two Gaussian moments
in the real case, or all Gaussian mixed moments of total degree at most
two in the complex case, and
$\sup_{N,i,j}\mathbb E|\xi_{ij}|^p\le M_p<\infty$.  Then
\begin{equation}\label{eq:moment-fractional-comparison}
 \left|\mathbb Em_X(i\eta)-\mathbb Em_{X^{\mathrm G}}(i\eta)\right|
 \le\frac{C_p(1+M_p)s_*^{(p-2)/2}}{\eta^{p+1}}.
\end{equation}
\end{proposition}

\begin{proof}
Use the notation of the previous proof and expand through degree two.
Let $f(x)=(2N)^{-1}\operatorname{Tr}R_x$, let $T_1(x),T_2(x)$ be the
degree-one and degree-two terms, and set
$\mathcal R_2(x)=f(x)-f(0)-T_1(x)-T_2(x)$.  The exact resolvent
remainder gives
\[
 |\mathcal R_2(x)|\le\frac{|x|^3}{N\eta^4}.
\]
Also $R_x-R_0$ has rank at most two and norm at most $2/\eta$.
Consequently
\[
 |f(x)-f(0)|\le\frac{2}{N\eta},\qquad
 |T_1(x)|\le\frac{|x|}{N\eta^2},\qquad
 |T_2(x)|\le\frac{|x|^2}{N\eta^3}.
\]
For $|x|\ge\eta$ this yields
$|\mathcal R_2(x)|\le4|x|^2/(N\eta^3)$.  Combining the two ranges gives
\begin{equation}\label{eq:moment-interpolated-remainder}
 |\mathcal R_2(x)|
 \le\frac4{N\eta}
       \min\left\{\left(\frac{|x|}{\eta}\right)^2,
                   \left(\frac{|x|}{\eta}\right)^3\right\}
 \le\frac{4|x|^p}{N\eta^{p+1}}.
\end{equation}
The degree-one and degree-two expectations cancel.  Summing the
remainder bound over entries, using
$N^{-1}\sum_{i,j}s_{ij}^{p/2}\le s_*^{(p-2)/2}$, proves the claim.
\end{proof}

\subsection{Coarse counts and fixed-scale convergence}

In the rest of the section assume the profile satisfies
Assumption~\ref{ass:profile-balance}.  For compactness, let $\theta>2$
denote either an integer $r\ge3$ satisfying $\mathcal M_r$, or an exponent
$p\in(2,3]$ satisfying Proposition~\ref{prop:moment-fractional}.  Define
\begin{equation}\label{eq:moment-exponents}
 \alpha_\theta:=\frac{\theta-2}{2(\theta+1)},\qquad
 q_\theta:=\frac1{1+\alpha_\theta}
          =\frac{2(\theta+1)}{3\theta}.
\end{equation}

\begin{corollary}[Non-Gaussian coarse count]\label{cor:moment-count}
There is a constant $K_\theta$, depending on the moment bound, $\theta$,
and $C_0$, such that for every fixed $z\ne0$,
\begin{equation}\label{eq:moment-count}
 \mathbb E F_{N,z}(t)\le C_{z,\theta}t,
 \qquad t\ge\eta_{0,\theta}:=K_\theta W^{-\alpha_\theta}.
\end{equation}
In particular,
\[
 \mathbb P\{F_{N,z}(t)>u t\}\le C_{z,\theta}/u
 \qquad(u>0, t\ge\eta_{0,\theta}).
\]
These are expectation and Markov estimates, not a fine-scale local law.
\end{corollary}

\begin{proof}
Since $\alpha_\theta<1/2$, choosing $K_\theta$ sufficiently large makes
every $t\ge\eta_{0,\theta}$ admissible for the Gaussian coarse estimate
in Theorem~\ref{thm:coarse}.  The relevant comparison proposition gives
\[
 \left|\mathbb Em_X(it)-\mathbb Em_{X^{\mathrm G}}(it)\right|
 \le C_\theta W^{-(\theta-2)/2}t^{-(\theta+1)}\le C_\theta,
\]
where constants may depend on the uniform moment bound and $C_0$.
The last inequality uses
$\alpha_\theta(\theta+1)=(\theta-2)/2$.  As
\[
 \operatorname{Im}m_X(it)=\frac1N\sum_{j=1}^N
          \frac{t}{s_j(X-zI_N)^2+t^2},
 \qquad F_{N,z}(t)\le2t\operatorname{Im}m_X(it),
\]
Theorem~\ref{thm:coarse} proves \eqref{eq:moment-count}.  Markov's
inequality proves the final assertion.
\end{proof}

\begin{proposition}[Non-Gaussian global singular-value law]
\label{prop:moment-global}
Under the preceding moment and profile assumptions, suppose $W\to\infty$.
For every fixed $z\in\mathbb C$, the measures $\nu_{N,z}$ converge weakly
in probability to the same deterministic measure $\nu_z$ as the Gaussian
reference in Proposition~\ref{prop:global}.
\end{proposition}

\begin{proof}
At every fixed $\eta>0$, the expectation comparison tends to zero because
$s_*\to0$ and $\theta>2$.  We establish concentration using
Efron--Stein rather than Gaussian Poincar\'e.  Resample a single entry
$x_{ij}$ independently, writing $\delta=x_{ij}'-x_{ij}$.  The
Hermitizations differ by $V(\delta)$, of rank at most two and norm
$|\delta|$.  The resolvent identity gives
\[
 |m_X(i\eta)-m_{X'}(i\eta)|
 \le\frac{|\delta|}{N\eta^2}.
\]
For complex-valued functions, Efron--Stein follows by applying its real
form to the real and imaginary parts.  Since
$\mathbb E|x_{ij}'-x_{ij}|^2=2s_{ij}$ and $\sum_{i,j}s_{ij}=N$,
\begin{equation}\label{eq:moment-efron-stein}
 \operatorname{Var}(m_X(i\eta))
 \le\frac12\sum_{i,j}\mathbb E
       |m_X(i\eta)-m_{X^{(ij)}}(i\eta)|^2
 \le\frac1{N\eta^4}.
\end{equation}
Here $X^{(ij)}$ means $X$ with only that entry independently resampled.
Thus for every rational $\eta\in(2,3)$ the integrals of
$h_\eta(s)=\eta/(s^2+\eta^2)$ converge in probability to their integrals
against $\nu_z$, by Proposition~\ref{prop:global} and the expectation
comparison.  Moreover,
\[
 \mathbb E\int s^2\,d\nu_{N,z}(s)=1+|z|^2.
\]
Exactly as in the tightness argument of Proposition~\ref{prop:global},
this makes the laws of the random measures tight.  Every subsequential
limit satisfies the same countable resolvent identities almost surely.
The Stieltjes-transform uniqueness argument given there forces that limit
to equal $\nu_z$.  Hence convergence is weakly in probability.
\end{proof}

\subsection{Circular-law consequences and their hypotheses}

\begin{theorem}[Moment hierarchy with a separate LSV input]
\label{thm:moment-circular}
Let the moment assumptions corresponding to $\theta$ and
Assumption~\ref{ass:profile-balance} hold, and suppose $W\to\infty$.
For almost every fixed $z\in\mathbb C\setminus\{0\}$, suppose there are
deterministic $L_N(z)\ge0$ and events $\Omega_{N,z}$ such that
\begin{equation}\label{eq:moment-LSV-hypothesis}
 \mathbb P(\Omega_{N,z})\to1,\qquad
 s_{\min}(X-zI_N)\ge e^{-L_N(z)}\text{ on }\Omega_{N,z},\qquad
 W^{-\alpha_\theta}L_N(z)\to0.
\end{equation}
Then $\mu_X$ converges weakly in probability to the circular law.
In particular, it is enough that, for every sufficiently small fixed
$\kappa>0$, such events exist with
\begin{equation}\label{eq:moment-LSV-power}
 L_N(z)\le C_zN^{1+\kappa}/W+O_z(1),
\end{equation}
and that $W\ge N^{q_\theta+c}$ for a fixed
$0<c<1-q_\theta$.
\end{theorem}

\begin{proof}
Corollary~\ref{cor:moment-count},
\eqref{eq:moment-LSV-hypothesis}, and Lemma~\ref{lem:log-ui} prove
lower logarithmic uniform integrability at the cutoff
$\eta_{0,\theta}=K_\theta W^{-\alpha_\theta}$.
The second-moment identity used in Proposition~\ref{prop:moment-global}
controls the upper logarithmic tail exactly as in
Section~\ref{sec:circular-proof}.  Proposition~\ref{prop:moment-global}
therefore gives convergence of the logarithmic potentials to the same
limit as normalized Ginibre.  The replacement principle from
Section~\ref{sec:circular-proof} proves the circular law.

For the final assertion, choose $0<\kappa<c(1+\alpha_\theta)$.  Since
$q_\theta(1+\alpha_\theta)=1$, \eqref{eq:moment-LSV-power} gives
\[
 W^{-\alpha_\theta}L_N(z)
 \le C_z\frac{N^{1+\kappa}}{W^{1+\alpha_\theta}}
       +O_z(W^{-\alpha_\theta})
 \le C_zN^{\kappa-c(1+\alpha_\theta)}+o(1)\longrightarrow0.
\]
This verifies \eqref{eq:moment-LSV-hypothesis}.
\end{proof}

\subsection{Density inputs and proof of the main moment theorem}\label{sec:standalone-lsv}

We use the Hilbert--Schmidt-truncated estimate of
\cite[Theorem~3.1, proved in Appendix~B]{han2026optimal}.
We record only the specialization needed here and its consequence;
the density and block-net argument is not repeated.

\begin{proposition}[Imported bounded-density estimate]
\label{newleastsingularvaluebound}
Let $X$ satisfy Definition~\ref{mainassumptionbands} and suppose
\begin{equation}\label{lowerboundsonvariances}
 s_{ij}\ge c_W/W\quad\text{whenever }|i-j|_N\le W
\end{equation}
for a fixed $c_W>0$. Assume the normalized atoms are real with uniformly
bounded densities, or complex with independent real and imaginary parts,
each having a uniformly bounded density. If $W\ge N^{1/2+c}$ for a
fixed $c>0$, then for every fixed $K,R>0$ and $0<\kappa<c/4$,
\begin{equation}\label{eq:standalone-lsv}
 \mathbb P\left\{s_{\min}(X-zI_N)
 \le t e^{-N^{3\kappa}N/W},\ \|X\|_{\mathrm{HS}}\le R\sqrt N\right\}
 \le C_{K,R,\kappa}t+e^{-N^{1+\kappa/4}},
 \qquad |z|\le K,\quad t>0,
\end{equation}
for all sufficiently large $N$. The constants may also depend on the
profile and density bounds and on $c$.
\end{proposition}

This is \cite[Theorem~3.1]{han2026optimal}. In the complex case, the
product of the two marginal density bounds gives the planar density
bound required there. We retain the density assumptions used in this
paper, although the cited theorem allows a broader complex class.
No symmetry of the profile or additional moment is required. The shift
is any deterministic complex number in the stated disk, including zero;
the uniformity does not assert simultaneous invertibility throughout a
disk for one realization. A lower band of width $c_*W$ and height
$c_*/W$ reduces to this formulation by a constant rescaling of $W$
and a slight decrease of the bandwidth slack.

\begin{corollary}[Removing the Hilbert--Schmidt cutoff]
\label{cor:polynomial-norm-lsv}
Under the assumptions of Proposition~\ref{newleastsingularvaluebound},
for every fixed $z\in\mathbb C$ and $0<\kappa<c/4$,
\begin{equation}\label{eq:second-moment-lsv-floor}
 \mathbb P\left\{s_{\min}(X-zI_N)
 \ge e^{-N^{3\kappa}N/W}\right\}\longrightarrow1.
\end{equation}
Only the normalized second moments are needed for this implication.
\end{corollary}

\begin{proof}
Apply \eqref{eq:standalone-lsv} with $\kappa/2$ and
$t_N=\exp\{-(N^{3\kappa}-N^{3\kappa/2})N/W\}$.
Since $W\le C_0N$, we have $t_N\to0$. For each fixed $R>0$,
normalization and Markov's inequality give
\[
 \mathbb P\{\|X\|_{\mathrm{HS}}>R\sqrt N\}
 \le\frac{\mathbb E\|X\|_{\mathrm{HS}}^2}{R^2N}=R^{-2}.
\]
Consequently the limsup, as $N\to\infty$, of the probability of
$s_{\min}(X-zI_N)<e^{-N^{3\kappa}N/W}$ is at most $R^{-2}$.
Letting $R\to\infty$ proves the assertion. The order of these limits
keeps $R$ fixed in the imported estimate and requires no uniform
integrability assumption on the family of normalized atoms.
\end{proof}
% END sections/lsv
% BEGIN sections/nongaussian

\begin{proof}[Proof of Theorem~\ref{thm:moment-main}]
In case \textup{(i)}, Proposition~\ref{prop:external-lsv} gives
\eqref{eq:moment-LSV-power} for every sufficiently small $\kappa>0$.
Theorem~\ref{thm:moment-circular} applies directly.
In case \textup{(ii)}, set $d=q_r+c-1/2>0$, so that
$W\ge N^{1/2+d}$.  Choose a fixed $\kappa>0$ sufficiently small that
\[
 \kappa<d/4,\qquad 3\kappa<c(1+\alpha_r).
\]
The second-moment floor \eqref{eq:second-moment-lsv-floor} in
Corollary~\ref{cor:polynomial-norm-lsv} supplies, for every fixed $z$,
events of probability tending to one on which
\[
 -\log s_{\min}(X-zI_N)\le L_N:=N^{3\kappa}N/W.
\]
The density and lower-band assumptions required there are precisely
those in case \textup{(ii)}; the Hilbert--Schmidt cutoff has already been
removed using the second moment.  Moreover,
\[
 W^{-\alpha_r}L_N
 =\frac{N^{1+3\kappa}}{W^{1+\alpha_r}}
 \le N^{3\kappa-c(1+\alpha_r)}\longrightarrow0.
\]
Theorem~\ref{thm:moment-circular} now proves the assertion.  No additional
moment assumption is used in this least-singular-value step.
\end{proof}

\begin{corollary}[Finite moment above two with a density input]
\label{cor:moment-finite-p}
Let $2<p\le3$ and suppose the hypotheses of
Proposition~\ref{prop:moment-fractional} and
Assumption~\ref{ass:profile-balance} hold.  Suppose also that either of
the two least-singular-value alternatives in
Theorem~\ref{thm:moment-main} holds.  Then the circular law holds whenever
\[
 W\ge N^{q_p+c},\qquad q_p=\frac{2(p+1)}{3p},\qquad
 0<c<1-q_p.
\]
Under the bounded-density and diagonal-lower-band alternative, the
uniform absolute $p$th moment is the only moment assumption beyond
centering and covariance normalization.
\end{corollary}

\begin{proof}
For the subgaussian alternative use
Proposition~\ref{prop:external-lsv}.  For the diagonal-lower-band
alternative use Corollary~\ref{cor:polynomial-norm-lsv}, choosing
$0<\kappa<(q_p+c-1/2)/4$ and
$3\kappa<c(1+\alpha_p)$.  In either case the same exponent calculation
as above verifies \eqref{eq:moment-LSV-hypothesis}, and
Theorem~\ref{thm:moment-circular} applies with $\theta=p$.
\end{proof}

For real entries, the first stages of this sufficient hierarchy are
\begin{center}
\begin{tabular}{c|c|c|c}
$r$ & extra real moments & $\alpha_r$ & $q_r$\\ \hline
$3$ & none & $1/8$ & $8/9$\\
$4$ & $\mathbb E\xi^3=0$ & $1/5$ & $5/6$\\
$5$ & $\mathbb E\xi^3=0$, $\mathbb E\xi^4=3$ & $1/4$ & $4/5$\\
$6$ & the preceding, and $\mathbb E\xi^5=0$ & $2/7$ & $7/9$
\end{tabular}
\end{center}
Each row also requires a uniformly bounded absolute $r$th moment for the
resolvent comparison, and an applicable least-singular-value estimate
for the circular-law conclusion.  For complex entries use mixed-moment
matching, not merely the displayed real identities.  With only a finite
$p$th moment, $2<p\le3$, the conditional exponent is
$q_p=2(p+1)/(3p)$, which tends to $1$ as $p\downarrow2$.

\begin{remark}[Scope of the hierarchy]
These are sufficient bandwidth exponents obtained from a finite
resolvent comparison, not claims of optimality.  Although
$q_r\downarrow2/3$ as $r\to\infty$, the proof is for each fixed $r$;
it does not justify taking $r$ to grow with $N$, nor does it prove the
Gaussian range for every $c>0$ for one fixed non-Gaussian law with only
finitely many matched moments.  The finite-$(2+\delta)$ comparison is likewise
not a least-singular-value theorem.  If neither the subgaussian-density
input nor the bounded-density diagonal-band input is available, one
must separately verify \eqref{eq:moment-LSV-hypothesis} before asserting
the circular law.
\end{remark}

% END sections/nongaussian
% BEGIN discrete-strip proofs
\section{Scalar strips: the elimination algorithm and its consequences}
\label{blsv:sec-tree}

We prove Theorems~\ref{blsv:thm-strip} and~\ref{blsv:thm-circle} in this
section, using one deformed-iid estimate whose proof is deferred to
Section~\ref{blsv:sec-local}. The algorithm is a recursive Schur-complement
elimination. Its purpose is to preserve independent square blocks until
the moment when an invertibility estimate is applied to them.

\subsection{What is retained, and why the flat core is enough}
\label{blsv:subsec-algorithm}
Write $\mathcal M=\sqrt W(X_N+B_N-zI)$. Two different geometric scales
are relevant. The support radius $W$ determines how wide a separator
must be to disconnect two intervals. The flat-core radius
$\lfloor W/q\rfloor$ determines the size of a square on which all entries
have the same law. A separator of size $W$ can therefore be covered by
$q$ consecutive diagonal squares of size comparable to $W/q$.
Inside each square the random part of $\mathcal M$ is
$\sqrt{a_N}G_j$, with $G_j$ iid. The entries outside these squares need
not have the same variance.

Here is one elimination step. Order a parent interval as its left child,
right child, and separator, denoted by $L,R,J$. Finite range gives
\[
 \mathcal M_{L\cup R\cup J}=
 \begin{pmatrix}
 A_L&0&U_L\\
 0&A_R&U_R\\
 V_L&V_R&D_J+\operatorname{diag}(\sqrt{a_N}G_1,\ldots,
                                  \sqrt{a_N}G_q)
 \end{pmatrix}.
\]
The matrix is not Hermitian: there is no relation between $U_L$ and
$V_L$, or between $U_R$ and $V_R$. First solve the two child problems.
Next reveal the coupling blocks and the entries included in $D_J$,
leaving the $G_j$ unrevealed. If both children are invertible, the
remaining matrix is
\[
 S_J=\operatorname{diag}(\sqrt{a_N}G_1,\ldots,\sqrt{a_N}G_q)
       +D_J-V_LA_L^{-1}U_L-V_RA_R^{-1}U_R.
\]
Conditional on everything revealed so far, the last three terms are
deterministic and the $G_j$ are still independent iid squares. The
Schur correction may be dense; no band structure is required of this
conditional deformation. This is the point at which the local
invertibility estimate applies.

For $q=1$, the entire separator is one fresh iid square. For $q=2$,
only the two diagonal half-squares are retained; all entries between
them are revealed. One first eliminates the first square and then
applies the same iid estimate to the second. The off-diagonal couplings
enter the conditional deformation; the local estimate imposes no
restriction on their support. It therefore applies to triangular or
incomplete couplings as well.
For Bernoulli atoms these retained entries are independent signs;
conditioning on the exterior does not change their law.

The recursion is executed from leaves to root:
\begin{enumerate}[label=\textup{\arabic*.}]
\item Remove a length-$W$ interval to open the circle. Bisect each long
remaining interval by another length-$W$ separator, stopping at size
at most $3W$.
\item At each leaf, reveal everything except a fixed number of diagonal
iid squares and eliminate those squares successively.
\item At a parent, use the solved children to form $S_J$, reveal the
retained squares successively, and recover the inverse of the parent
from the block inverse formula.
\end{enumerate}
The root has a single exterior interval; the same Schur formula applies
with one child term omitted. Balanced splitting makes the number of
successive parent steps logarithmic in $N/W$.

\subsection{The local estimate and the cost of one separator}
Throughout the discrete proof, $\xi$ is the fixed centered, unit-variance
subgaussian atom of Part II; it may be real or complex. The following
estimate is the only Littlewood--Offord input needed by the algorithm.
\begin{lemma}[Fixed-power quantitative invertibility]\label{blsv:lem-local}
Let $G_n$ have independent entries with law $\xi$, and let $D$ be any
deterministic complex $n\times n$ matrix with $\|D\|\le e^{n^{1/5}}$.
Set $H=\|D\|+C_\xi\sqrt n$, increasing $C_\xi$ so that $H\ge\sqrt n$.
There is $C<\infty$ such that
\begin{equation}\label{blsv:eq-local-result}
 \blsvP\{\blsvsmin(D+G_n)\le(CHn^3)^{-3}\}\le C n^{-1/8}.
\end{equation}
\end{lemma}

Its small-ball and counting ingredients are those of Jain. The
reparameterization gives a cubic deformation loss with failure
probability $O(n^{-1/8})$. Its proof in Section~\ref{blsv:sec-local}
can be read independently of the strip geometry.

We make the deterministic norm bookkeeping explicit. If
$M=\left(\begin{smallmatrix}A&B\\ C&D\end{smallmatrix}\right)$,
$A$ is invertible, and $S=D-CA^{-1}B$ is invertible, then the block inverse
formula gives
\begin{equation}\label{blsv:eq-block-inverse}
 \|M^{-1}\|\le \|A^{-1}\|+
 (1+\|A^{-1}B\|)(1+\|CA^{-1}\|)\|S^{-1}\|.
\end{equation}
For example, the second term is the norm bound for
$\binom{-A^{-1}B}{I}S^{-1}(-CA^{-1},I)$.

For orientation, suppose a child inverse is bounded by $R$ and all
couplings have polynomial norm in $N$. When $q=1$, the Schur deformation
has norm at most $N^C R$, the local lemma gives a Schur inverse at most
$N^C R^3$, and the two factors in \eqref{blsv:eq-block-inverse} cost
another $R^2$. Thus the parent bound is $N^C R^5$. With $q$ fresh
squares the same bookkeeping gives $N^C R^{5^q}$; the next lemma proves
the finite-dimensional step precisely.

\begin{lemma}[A fixed number of fresh iid squares]\label{blsv:lem-finite}
Fix an integer $r\ge1$, $b>0$, and $0<c<C<\infty$. Let
$bW\le n_j\le b^{-1}W$ for $1\le j\le r$, and let $G_j$ be mutually
independent iid matrices of sizes $n_j$, with atom $\xi$. Let
$c\le|\sigma_j|\le C$, and set
\[
 M=D+\operatorname{diag}(\sigma_1G_1,\ldots,\sigma_rG_r),
 \qquad \|D\|\le H,\quad H\ge2.
\]
There are $c_r,C_r>0$ such that, if
\begin{equation}\label{blsv:eq-finite-cap}
 \log H+\log W\le c_r W^{1/5},
\end{equation}
then, except on an event of probability at most $C_rW^{-1/8}$,
\begin{equation}\label{blsv:eq-finite-inverse}
 M\text{ is invertible},\qquad
 \|M^{-1}\|\le W^{C_r}H^{5^r-2}.
\end{equation}
The estimates are uniform over deterministic complex $D$.
\end{lemma}

\begin{proof}
For $r=1$, divide by $\sigma_1$ and apply
Lemma~\ref{blsv:lem-local}. Its bound is at most $W^{C_1}H^3$;
the $\sqrt{n_1}$ term is absorbed into the polynomial in $W$.

More explicitly, proceed inductively through the leading principal
blocks. Suppose the first $j$ blocks have inverse norm
$R\le W^{A_j}H^{a_j}$. Their coupling to the next block has norm at most
$H$. Conditional on these blocks, the next Schur complement is
$\sigma_{j+1}G_{j+1}+D'$ with
\[
 \|D'\|\le H+H^2R\le W^{A'_j}H^{a_j+2}.
\]
The deformation $D'$ is measurable with respect to the preceding
blocks and is independent of the fresh $G_{j+1}$.
Choose $c_r$ sufficiently small in \eqref{blsv:eq-finite-cap}.
For every $j\le r$, the displayed deformation after division by
$\sigma_{j+1}$ has norm at most $e^{n_{j+1}^{1/5}}$.
Lemma~\ref{blsv:lem-local} therefore gives Schur inverse norm at most
$W^{A''_j}H^{3a_j+6}$ with conditional failure $CW^{-1/8}$.
By \eqref{blsv:eq-block-inverse}, the new inverse norm is bounded by
$W^{A_{j+1}}H^{5a_j+8}$. The exponents satisfy
\[
 a_1=3,\qquad a_{j+1}=5a_j+8,\qquad a_j=5^j-2.
\]
There are only $r$ steps; a conditional union bound proves the assertion.
\end{proof}

\subsection{Proof of the strip estimate: geometry, exposure, and recursion}
\begin{proof}[Proof of Theorem~\ref{blsv:thm-strip}]
Work with
\[
 \mathcal M=\sqrt W\,(X_N+B_N-z\blsvId_N),
\]
where $B_N=0$ is allowed. The constants below may depend on
$q,c_0,C_0,K,A,\xi$. Since the atom is subgaussian, except on an event
of probability at most $CN^2e^{-cN^2}$, every active $\xi_{i,s}$ has
absolute value at most $N$. On this event all off-diagonal coupling
blocks of $\mathcal M$ have norm at most $N^{A_0}$ for a fixed $A_0$.
We refer to these estimates as the polynomial coupling bounds. At each
conditioning step below, we impose only the bounds on already exposed
couplings, so the unexposed iid blocks retain their original laws.

\smallskip\noindent\emph{Geometry of the tree.}
Remove one consecutive interval of length $W$ from the circle. It is
the root separator. The complement has the support of a linear
finite-range strip: a connection through the removed interval would
have length at least $W+1$. For a remaining interval of length $m>3W$,
remove a central interval of length $W$. Its two children have lengths
$\lfloor(m-W)/2\rfloor$ and $\lceil(m-W)/2\rceil$, both at least $W$.
There is no matrix entry connecting the two children. Stop when an
interval has length at most $3W$.

The tree has at most $CN/W$ nodes and height at most
\begin{equation}\label{blsv:eq-height}
 h\le\log_2(N/W)+C.
\end{equation}
These assertions follow either from the balanced splitting or by
charging each disjoint separator or leaf to at least $W$ coordinates.

Each length-$W$ separator can be partitioned into $q$ consecutive
intervals with sizes differing by at most one. Their maximum diameter
is $\lceil W/q\rceil-1\le\lfloor W/q\rfloor$.
Consequently the random matrix on each diagonal square is
$\sqrt{a_N}G_j$, a full iid square. All remaining separator entries may
be exposed. These $q$ squares are independent of one another and of all
entries used to construct the exterior Schur correction.

Similarly, any leaf of size between $W$ and $3W$ can be partitioned into
at most $3q+1$ consecutive intervals, each of size comparable to $W$
with constants depending on $q$, and diameter at most
$\lfloor W/q\rfloor$. For example, partition it as evenly as possible
into $\lceil m/(\lfloor W/q\rfloor+1)\rceil$ intervals.

\smallskip\noindent\emph{Initialization at the leaves.}
Expose the off-diagonal blocks of this leaf partition. The fresh
diagonal squares have the form required by Lemma~\ref{blsv:lem-finite},
and the deterministic deformation has polynomial norm in $N$ on the
polynomial coupling event. Condition \eqref{blsv:eq-cap-global} implies
$\log N=o(W^{1/5})$. Thus that lemma, with at most $3q+1$ blocks, gives
leaf inverse norm at most $N^{A_1}$, with failure at most $CW^{-1/8}$.

\smallskip\noindent\emph{One parent step.}
Suppose its children are invertible with inverse norms at most $R\ge1$.
The matrix on their union is block diagonal, so its inverse norm is at
most $R$. After eliminating it, the Schur complement on the separator
has deterministic deformation norm at most
\begin{equation}\label{blsv:eq-parent-H}
 H\le N^{A_2}R,
\end{equation}
once the inter-square entries in the separator are exposed.
The $q$ diagonal iid squares remain fresh. Lemma~\ref{blsv:lem-finite}
then bounds its inverse by $N^{A_3}R^{5^q-2}$, with conditional failure
at most $CW^{-1/8}$. Applying \eqref{blsv:eq-block-inverse} to the
exterior elimination gives
\begin{equation}\label{blsv:eq-parent-recurrence}
 R_{\rm parent}\le N^{A_4}R^{5^q}.
\end{equation}

\smallskip\noindent\emph{Checking every deformation cap.}
Set $\lambda=5^q$ and start from $R_0=N^{A_1}$. The deterministic
recursion $R_{j+1}=N^{A_4}R_j^\lambda$ satisfies
\[
 \log R_j\le C\lambda^j\log N.
\]
By \eqref{blsv:eq-height}, every $\log R_j$ and the logarithm of the
upper bound in \eqref{blsv:eq-parent-H} are at most
$C\Theta_{N,W,q}=o(W^{1/5})$. Therefore
\eqref{blsv:eq-finite-cap} holds at every node for all sufficiently
large $N$. This also verifies the caps used at the leaves; the number
of subblocks at a leaf is fixed, so its larger constant is harmless.

\smallskip\noindent\emph{Probability bookkeeping.}
At an internal node, condition on all entries outside its fresh diagonal
squares. The child success events and the required coupling bounds are
measurable under this conditioning. The conditional probability that
the node fails its asserted inverse bound while these prerequisites
hold is at most $CW^{-1/8}$. The same bound holds for a leaf.
If all the coupling bounds hold but the recursive conclusion fails,
there is a first failed node with successful prerequisites.
Thus a union bound over at most $CN/W$ nodes costs at most
$CN/W^{9/8}$. Adding the discarded coupling event does not change this
bound for large $N$.

On the resulting success event,
$\|\mathcal M^{-1}\|\le\exp(C\Theta_{N,W,q})$.
Since $X_N+B_N-zI=W^{-1/2}\mathcal M$, the additional factor $\sqrt W$
in its inverse is absorbed by increasing $C$. A further enlargement
of $C$ makes the success lower bound strictly larger than the threshold
in \eqref{blsv:eq-lsv-global}, so that the weak inequality in its bad
event causes no endpoint issue. Every estimate is uniform for
deterministic $|z|\le K$.
\end{proof}

\begin{remark}[Local scaling and limits of the profile extension]
The same elimination argument applies when each fresh square is
$P_jG_jQ_j$ with deterministic diagonal $P_j,Q_j$ whose norms and inverse
norms are bounded uniformly: divide on both sides before applying the
local lemma. It also permits arbitrary independent exterior couplings
with a polynomial norm bound when only invertibility is sought.
The argument does not prove the same bound for arbitrary independently
weighted entries inside every diagonal square. In particular, a
two-sided variance bound alone does not turn the square into an iid
matrix. Nor may a noncompact profile simply be truncated without a
separate estimate on the effect of that truncation at the very small
singular-value scale.
\end{remark}

\subsection{From the inverse bound to the Bernoulli circular law}
\label{blsv:sec-circle}

\begin{proof}[Proof of Theorem~\ref{blsv:thm-circle}]
The flat-core profile is doubly stochastic and invariant under cyclic
shifts, so Assumption~\ref{ass:profile-balance} holds even when
$p_{N,s}\ne p_{N,-s}$. For real subgaussian atoms, the comparison results
of Section~\ref{sec:nongaussian} apply with $r=3$, or with $r=4$ when
$\mathbb E\xi^3=0$. Their counting cutoff is $W^{-\alpha}$, where
$\alpha=1/8$ or $1/5$, respectively; those comparison results require
no density assumption.

Theorem~\ref{thm:moment-circular} reduces the proof to verifying
\eqref{eq:moment-LSV-hypothesis}. The strip estimate supplies
\[
 L_N=C(N/W)^{\nu_q}\log N,\qquad
 \mathbb P(\Omega_{N,z}^c)\le CN/W^{9/8}.
\]
Consequently the only conditions to check are
\begin{equation}\label{blsv:eq-log-condition}
 N/W^{9/8}\to0,\qquad
 W^{-\alpha}(N/W)^{\nu_q}\log N\to0,\qquad
 (N/W)^{\nu_q}\log N=o(W^{1/5}).
\end{equation}
The second condition is the decisive one: the coarse count permits
mass $O(W^{-\alpha})$ below its cutoff, while the inverse estimate
bounds the logarithmic cost per singular value by $L_N$.
Lemma~\ref{lem:log-ui} controls this product; the remaining logarithmic
tail and the replacement principle were already handled in
Theorem~\ref{thm:moment-circular}.

If $W\ge N^\gamma$ with $\gamma>\nu_q/(\nu_q+\alpha)$, then
\[
 W^{-\alpha}(N/W)^{\nu_q}\log N
 \le N^{\nu_q(1-\gamma)-\alpha\gamma}\log N\longrightarrow0.
\]
Since $\alpha\le1/5$, this also verifies the deformation cap, the last
condition in \eqref{blsv:eq-log-condition}. Moreover,
$\nu_q\ge\log_2 5$ and $\alpha\le1/5$ give
$\nu_q/(\nu_q+\alpha)>8/9$, so the failure probability tends to zero.
All hypotheses of Theorem~\ref{thm:moment-circular} are satisfied.
\end{proof}

For Rademacher atoms, $\mathbb E\xi^3=0$, so $\alpha=1/5$ applies.
The role of the two parts of the paper is now explicit: moment
comparison supplies the count of small singular values, and the
separator algorithm supplies their lowest admissible scale. The latter
uses the finite-range and flat-core assumptions; the argument does not
extend the density-free conclusion to the general graph profiles in
Part I.

% END discrete-strip proofs
\section{Proof of the quantitative deformed-iid estimate}\label{blsv:sec-local}

Section~\ref{blsv:sec-tree} uses only the statement of
Lemma~\ref{blsv:lem-local}. We now prove it. The aim is a fixed cubic
loss in the deformation norm, since that power enters the separator
recursion. We use Jain's existing inverse Littlewood--Offord inputs;
the change is in the choice of scales and parameters.

The proof separates unit vectors according to the concentration of
their random linear forms. Poor vectors are handled by conditioning on
all but one row. Rich vectors are rounded to integer vectors; Jain's
counting theorem makes the resulting set small enough for a union bound
over independent rows. Three scales ensure that at least one rounding
scale has controlled growth of its concentration function.

We record precisely the inputs from Jain~\cite{blsvJain}. For
$v\in\blsvC^n$ and $t\ge0$, write
\[
 \blsvQ_t(v)=\sup_{b\in\blsvC}
    \blsvP\left\{\left|\sum_{j=1}^n v_j\xi_j-b\right|\le t\right\}.
\]
Multiplication of $v$ by a complex scalar $d\ne0$ gives
$\blsvQ_t(dv)=\blsvQ_{t/|d|}(v)$.

\subsection{Imported anti-concentration and counting inputs}

Every fixed nondegenerate unit-variance atom is $C_\xi$-good in the sense
of Jain's Definition~2.6. We use the following statements, with constants
depending only on the atom.

\begin{enumerate}
\item By Jain's Lemma~2.3, there is $c_*>0$ such that
\begin{equation}\label{blsv:eq-uniform-ac}
 \sup_{\|v\|=1}\blsvQ_{c_*}(v)\le1-c_*.
\end{equation}
\item Jain's Proposition~2.8 implies the following rounding statement.
If $\blsvQ_1(w)>\beta$, $f=c\beta$, and $g=A=n^\delta$ for any fixed
$\delta>0$, then, for all large $n$ and sufficiently small $c=c_\xi>0$,
there are $d\in\blsvC$ and $u\in(\blsvZ+i\blsvZ)^n$ such that
\begin{equation}\label{blsv:eq-round-input}
 f\le|d|\le n^\delta,\qquad \|dw-u\|\le n^\delta,
\end{equation}
provided $\beta$ is a fixed inverse power of $n$. Indeed, the contrary
would bound $\blsvQ_1(w)$ by
$C e^\pi(2e^{-c n^{2\delta}}+f)<\beta$.
\item For
$\mathcal V_\rho=\{u\in(\blsvZ+i\blsvZ)^n:\blsvQ_1(u)\ge\rho\}$,
Jain's Theorem~1.3 gives
\begin{equation}\label{blsv:eq-jain-count}
 |\varphi_p(\mathcal V_\rho)|
 \le\left(\frac{5np^2}{s}\right)^s+
 \left(\frac{C\rho^{-1}}{\sqrt{s/k}}\right)^n,
\end{equation}
if $1000C_\xi\le k\le\sqrt{s}\le s\le n/\log n$, the odd prime $p$
satisfies $C\rho^{-1}\le p\le2^{n/s}$, and
$\rho\ge C\max\{e^{-s/k},s^{-k/4}\}$.
\end{enumerate}

Two elementary subgaussian estimates will also be used. For deterministic
$v,w$,
\begin{equation}\label{blsv:eq-perturb-q}
 \blsvQ_t(v)\le\blsvQ_{t+r}(w)
       +C\exp\{-cr^2/\|v-w\|^2\},\qquad r>0,
\end{equation}
with the error term omitted when $v=w$. This follows by separating the
event $|\sum(v_j-w_j)\xi_j|>r$. Furthermore, an iid matrix $G_n$ with
atom $\xi$ satisfies
\begin{equation}\label{blsv:eq-norm-iid}
 \blsvP\{\|G_n\|>C_\xi\sqrt n\}\le C_\xi e^{-c_\xi n}.
\end{equation}
These are also the inputs in Jain's subgaussian proof. Covering a disk
of radius $R\ge1$ by $CR^2$ disks of radius $1$ yields
\begin{equation}\label{blsv:eq-cover-q}
 \blsvQ_1(v)\ge cR^{-2}\blsvQ_R(v).
\end{equation}

\subsection{Three-scale rounding and the counting budget}
The numerical parameters below serve two purposes. The poor-vector
union bound costs $n\beta=n^{-1/8}$. For rich vectors, three scales
force a concentration plateau with loss $n^{0.38}$, while counting
produces a gain $n^{-0.384}$ per coordinate. The remaining factor
$n^{-0.004n}$ pays for rounding and exceptional rows. The smaller
auxiliary exponents leave room for the finite-field counting conditions;
they do not represent additional geometric restrictions on the strip.

\begin{proof}[Proof of Lemma~\ref{blsv:lem-local}]
Write $F=D+G_n$ and fix
\begin{equation}\label{blsv:eq-local-parameters}
 \eta=(C_1Hn^3)^{-3},\quad \beta=n^{-9/8},\quad f=c_1\beta,
 \quad\Lambda=C_2H/f,\quad r_j=2\eta\sqrt n\,\Lambda^j\ (0\le j\le3).
\end{equation}
Here $c_1$ is sufficiently small, $C_2$ is sufficiently large, and $C_1$
is then chosen sufficiently large. In particular,
\begin{equation}\label{blsv:eq-r3}
 r_3\le C n^{-41/8},\qquad r_3=o(1).
\end{equation}
All constants below are independent of $D$ and $n$.

\smallskip\noindent\emph{Poor vectors.}
Call a unit vector poor if $\blsvQ_{r_0}(v)\le\beta$.
We give the row-conditioning argument, so no restriction from the
numerical statement of Jain's main theorem is being imported here.
For each row index $i$, condition on all other rows. If there is a poor
unit vector $v$ with $\|F_{-i}v\|\le\eta$, select one, denoted $v^{(i)}$,
using a measurable selection on this feasible set. The feasible set is
Borel in the conditioning data and the vector; universal measurable
selection suffices on the completed probability space. If it is empty,
the event for this index is discarded.

Suppose a poor unit vector $v$ satisfies $\|Fv\|\le\eta$. Then there is
a unit left singular vector $w$ with $\|w^*F\|\le\eta$. Choose $i$ with
$|w_i|\ge n^{-1/2}$. The feasible set for $i$ is nonempty, and
\[
 |\overline{w_i}F_i v^{(i)}|
 \le |w^*Fv^{(i)}|+
       \left|\sum_{j\ne i}\overline{w_j}F_jv^{(i)}\right|
 \le2\eta.
\]
Thus $|F_iv^{(i)}|\le r_0$. Conditional on the other rows, its probability
is at most $\beta$, since $v^{(i)}$ is poor and the deterministic part of
row $i$ only changes the center. A union bound gives
\begin{equation}\label{blsv:eq-poor}
 \blsvP\{\exists\text{ poor }v:\|Fv\|\le\eta\}\le n\beta=n^{-1/8}.
\end{equation}

\smallskip\noindent\emph{Three scales for rich vectors.}
For a rich unit vector, $\blsvQ_{r_0}(v)>\beta$. Since
$3(0.38)=1.14>9/8$, some $j\in\{0,1,2\}$ satisfies
\begin{equation}\label{blsv:eq-plateau}
 \blsvQ_{r_{j+1}}(v)\le n^{0.38}\blsvQ_{r_j}(v).
\end{equation}
Otherwise $\blsvQ_{r_3}(v)>n^{1.14}\beta>1$ for large $n$.
Partition the possible values $\blsvQ_{r_j}(v)$ into dyadic intervals
$[\vartheta,2\vartheta]$, with $\vartheta\ge\beta/2$; there are $O(\log n)$ such intervals.
The endpoints can be assigned to either adjacent interval.

Apply \eqref{blsv:eq-round-input} to $w=v/r_j$, with $\delta=0.001$.
There exist $d$ and a Gaussian integer vector $u$ with
\begin{equation}\label{blsv:eq-round}
 f\le|d|\le n^{0.001},\qquad
 v'=dv/r_j,\qquad \|v'-u\|\le n^{0.001}.
\end{equation}
For each fixed pair $(j,\vartheta)$ let $\mathcal U_{j,\vartheta}$ be the set of all
integer vectors obtainable in this way from its rich vectors and valid
choices of $d$. This is a deterministic finite set; it does not depend
on the realization of $G_n$.

\smallskip\noindent\emph{Counting the rounded candidates.}
Since $\blsvQ_{|d|}(v')=\blsvQ_{r_j}(v)\ge \vartheta$, the subgaussian
comparison with error radius $n^{0.002}$ gives
\[
 \blsvQ_{2n^{0.002}}(u)\ge \vartheta-Ce^{-c n^{0.002}}\ge \vartheta/2
\]
for large $n$. By \eqref{blsv:eq-cover-q},
\begin{equation}\label{blsv:eq-q1-u}
 \blsvQ_1(u)\ge c \vartheta n^{-0.004}.
\end{equation}
Moreover,
\[
 \|u\|_\infty\le n^{0.001}/r_0+n^{0.001},\qquad
 \log(2+\|u\|_\infty)\le 3\log H+C\log n=O(n^{0.2}).
\]
Take integers $s=\lfloor n^{0.78}\rfloor$, $k=\lfloor n^{0.004}\rfloor$
and an odd prime $p\in[e^{n^{0.21}},2e^{n^{0.21}}]$. For all large $n$
the map $\varphi_p$ is injective on all candidates: both the real and
imaginary coordinates have absolute value less than $p/3$.
All conditions of \eqref{blsv:eq-jain-count} hold with
$\rho=c\vartheta n^{-0.004}$. In particular $\rho\ge c n^{-9/8-0.004}$,
$k\le\sqrt{s}$, and $\log p=O(n^{0.21})=o(n/s)$.
It follows that
\begin{equation}\label{blsv:eq-candidate-count}
 |\mathcal U_{j,\vartheta}|
 \le \exp(Cn^{0.99})+
       (C\vartheta^{-1}n^{-0.384})^n.
\end{equation}
Here $0.384=(0.78-0.004)/2-0.004$, and $0.99=0.78+0.21$.

\smallskip\noindent\emph{One-row bounds for the same candidates.}
Fix $u\in\mathcal U_{j,\vartheta}$ with one witnessing pair $(v,d)$.
Using \eqref{blsv:eq-perturb-q} with error radius $H$, and enlarging $C_2$
in \eqref{blsv:eq-local-parameters} if necessary, yields
\begin{align}
 \blsvQ_{C_3H}(u)
 &\le \blsvQ_{(C_3+1)H}(v')+Ce^{-cH^2/n^{0.002}}\notag\\
 &\le \blsvQ_{r_{j+1}}(v)+Ce^{-c n^{0.998}}.
 \label{blsv:eq-row-comparison}
\end{align}
Here $C_3$ is any fixed constant needed below, and
$(C_3+1)Hr_j/|d|\le r_{j+1}$ follows from $|d|\ge f$.
Equations \eqref{blsv:eq-plateau}, \eqref{blsv:eq-r3}, and
\eqref{blsv:eq-uniform-ac} imply, for large $n$,
\begin{equation}\label{blsv:eq-row-bound}
 \blsvQ_{C_3H}(u)\le\min\{1-c,\ C\vartheta n^{0.38}\}.
\end{equation}
The exponentially small error is absorbed because $\vartheta\ge\beta/2$.
In particular no candidate used here is zero.

\smallskip\noindent\emph{Eliminating rich near-kernel vectors.}
On $\mathcal O=\{\|G_n\|\le C_\xi\sqrt n\}$, one has $\|F\|\le H$.
If a rich unit vector $v$ satisfies $\|Fv\|\le\eta$, its candidate $u$
satisfies
\[
 \|Fu\|\le |d|\eta/r_j+Hn^{0.001}\le C Hn^{0.001}.
\]
Consequently at least $n-\ell_n$ rows satisfy $|(Fu)_i|\le C_3H$,
where $\ell_n=\lceil Cn^{0.002}\rceil$.
For a deterministic candidate the rows are independent. Summing over
their possible exceptional sets costs at most
$\exp(C\ell_n\log n)$. Combining this with
\eqref{blsv:eq-candidate-count} and \eqref{blsv:eq-row-bound}, the first
candidate-count term contributes at most
\[
 \exp(Cn^{0.99}+C\ell_n\log n)(1-c)^{n-\ell_n}\le e^{-c'n}.
\]
For the second term use the bound $C\vartheta n^{0.38}$ in
\eqref{blsv:eq-row-bound}. Its contribution is at most
\[
 \exp(C\ell_n\log n)
 (C\vartheta^{-1}n^{-0.384})^n(C\vartheta n^{0.38})^{n-\ell_n}
 \le \exp\{-0.004n\log n+Cn+C\ell_n\log n\}.
\]
The last estimate uses $\vartheta\ge n^{-9/8}/2$. Summing over the three values
of $j$ and the $O(\log n)$ dyadic intervals, and adding
$\blsvP(\mathcal O^c)$, bounds the rich-vector contribution by $Ce^{-cn}$.
Together with \eqref{blsv:eq-poor}, this proves
\eqref{blsv:eq-local-result}.
\end{proof}

\appendix
\section{The former claims as open problems, and the failure of their proof}
\label{sec:former-claims}
\label{sec:conjectures}

In \href{https://arxiv.org/abs/2508.18143v2}{arXiv:2508.18143v2}, we claimed circular-law convergence at bandwidth
\(W\geq N^{1/2+c}\), together with a fine-scale entrywise local law.
The proof of that local law contains an incorrect sign in its linearized
self-consistency equation.  Correcting the sign changes the relevant
stability operator and does not complete the argument.  We therefore
withdraw the former proofs and retain the statements below as conjectures
in their original generality, not as results proved by that argument.
Specific subclasses of the circular-law statements are now covered by
\cite{han2025slow,han2026optimal}; the scope of those results was described
in the introduction. The fine local law itself remains an open problem.
Theorem~\ref{thm:circular-gaussian}
uses the coarse resolvent estimate proved separately and does not depend
on any conjecture in this appendix.

\subsection{Historical hypotheses and the full former statements}

All former theorem and definition numbers in this appendix refer to
arXiv:2508.18143v2. For clarity, we restate the hypotheses used in \href{https://arxiv.org/abs/2508.18143v2}{arXiv:2508.18143v2},
including its inverse condition.  Keeping this condition here records the
scope of the former assertions; it does not assert that it is the correct
stability condition for their proof.

\begin{assumption}[Historical variance-profile and moment hypotheses]
\label{ass:historical-profile}
Let \(X_{ij}=\sqrt{s_{ij}}\,\xi_{ij}\), where the variables \(\xi_{ij}\)
are independent, \(\mathbb E\xi_{ij}=0\), and
\(\mathbb E|\xi_{ij}|^2=1\).  The deterministic matrix
\(S=(s_{ij})_{i,j=1}^N\) is nonnegative and doubly stochastic, and
\[
 \max_{i,j}s_{ij}\leq \frac{C}{W}.
\]
For each \(z\) with \(0<|z|<1\), there are constants
\(c_{S,z}>0\) and \(C_{S,z}>0\), independent of \(N,W\), such that
\begin{equation}
 \left\|
 \left(I_{2N}-
 \begin{pmatrix}
  y_2S&y_1S^T\\
  y_1S&y_2S^T
 \end{pmatrix}\right)^{-1}
 \right\|_{\ell^\infty\to\ell^\infty}
 \leq C_{S,z}(\log N)^{C_{S,z}}
 \label{eq:historical-inverse}
\end{equation}
whenever
\[
 |y_1+1-|z|^2|\leq c_{S,z},\qquad
 |y_2+|z|^2|\leq c_{S,z}.
\]
The constants may depend on \(|z|\) and the fixed profile family.
Finally, for every fixed integer \(p\geq1\),
there is a constant \(\mu_p<\infty\), independent of \(N,i,j\), such that
\[
 (\mathbb E|\xi_{ij}|^p)^{1/p}\leq\mu_p.
\]
\end{assumption}

The former Definition~1.1 was not restricted to symmetric profiles.
Its ``bandwidth'' hypothesis was the upper bound \(s_{ij}\leq C/W\),
not a support condition.  Its constants in
\eqref{eq:historical-inverse} were allowed to depend on \(|z|\).
In particular, none of the following statements silently
replaces the old profile class by a flat full-matrix model.

We say that the entries satisfy the \emph{historical density condition}
if either all \(\xi_{ij}\) are real and have densities bounded uniformly
on \(\mathbb R\), or they are complex with independent real and imaginary
parts, each having a uniformly bounded density on \(\mathbb R\).
Uniform subgaussian tails mean that
\(\mathbb E\exp(|\xi_{ij}|^2/K^2)\leq2\) for a common constant \(K\).

\begin{conjecture}[Former circular-law claims]
\label{conj:historical-circular}
For every fixed \(c>0\), the empirical eigenvalue measure of \(X\)
converges in probability to the uniform probability measure on the unit
disk whenever \(W\geq N^{1/2+c}\), in each of the following settings.
\begin{enumerate}[label=\textup{(\roman*)}]
 \item The profile satisfies the historical conditions in
 Assumption~\ref{ass:historical-profile}, and the entries satisfy the
 historical density condition and have uniform subgaussian tails.
 This is the assertion of former Theorem~1.6.
 \item Assumption~\ref{ass:historical-profile} and the historical density
 condition hold, and there is a fixed \(c_*>0\) such that
 \[
  s_{ij}\geq\frac{c_*}{W}
  \quad\text{whenever}\quad |i-j|_N\leq c_*W.
 \]
 This is the circular-law assertion in former Theorem~1.8.
\end{enumerate}
\end{conjecture}

Only the circular-law part of former Theorem~1.8 is included in this
conjecture.  Its separate least-singular-value assertion is not a
consequence of the local-law argument discussed here.  Conversely, an
exponential lower bound on the least singular value does not by itself
prove any of the circular-law claims above.

Here is the full fine-scale local-law assertion of former Theorem~3.4.
We retain the Gram spectral parameter used there, to distinguish it
from the Hermitized resolvent \(R\) in the coarse estimate.
For \(Y_z=X-zI\), set
\begin{equation}
 G(w,z)=(Y_z^*Y_z-wI)^{-1},\qquad
 \mathcal G(w,z)=(Y_zY_z^*-wI)^{-1}.
 \label{eq:historical-Gram-resolvents}
\end{equation}
Let \(m_c(w,z)\) be the Stieltjes-transform solution of
\begin{equation}
 \frac1{m_c}=-w(1+m_c)+\frac{|z|^2}{1+m_c},
 \label{eq:historical-scalar-equation}
\end{equation}
with positive imaginary part and normalization \(m_c(w,z)\sim-1/w\)
at infinity in the upper half-plane.

\begin{conjecture}[Full former fine-scale local-law statement]
\label{conj:historical-local-law}
Suppose Assumption~\ref{ass:historical-profile} holds and
\(W\geq N^{c_0}\) for a fixed \(c_0>0\).  Fix a sufficiently small
\(\tau>0\) and \(z\) with \(\tau\leq|z|\leq1-\tau\).
For a fixed sufficiently small \(\gamma_0>0\), define
\begin{equation}
 \mathbf S_{\gamma_0}
 =\{\eta>0:N^{\gamma_0}W^{-2}\leq\eta\leq10\}.
 \label{eq:historical-eta-domain}
\end{equation}
For every sufficiently small \(\varepsilon>0\) and every \(D>0\),
for all sufficiently large \(N\),
\begin{equation}
 \mathbb P\left(
  \bigcap_{\eta\in\mathbf S_{\gamma_0}}
  \bigcap_{1\leq i,j\leq N}
  \left\{
   |G_{ij}(i\eta,z)-m_c(i\eta,z)\delta_{ij}|
   \leq N^\varepsilon W^{-1/2}\eta^{-3/4}
  \right\}
 \right)\geq1-N^{-D}.
 \label{eq:historical-local-law-conjecture}
\end{equation}
The constants may depend on the fixed parameters and moment bounds,
but not on \(N,W,i,j\).
\end{conjecture}

The range in \eqref{eq:historical-eta-domain} and the error in
\eqref{eq:historical-local-law-conjecture} are those of the previous
statement, not those of the proved coarse estimate.  There is an
ambiguity in the old bandwidth declarations that we make explicit:
former Theorem~3.4 referred to the hypotheses of former Theorem~2.1,
in a section assuming only \(W\geq N^{c_0}\), whereas its proof also
said ``Recall that \(W\geq N^{1/2}\).''  We have recorded the broader
declared statement as an open problem.  Even its restriction to
\(W\geq N^{1/2+c}\), the range needed for the former circular-law
claims, is not established by that proof.

\begin{conjecture}[Former intermediate singular-value count]
\label{conj:historical-count}
Under Assumption~\ref{ass:historical-profile}, with
\(W\geq N^{c_0}\) for any fixed \(c_0>0\), fix
\(\tau\leq|z|\leq1-\tau\).  For every \(D,\varepsilon>0\),
\begin{equation}
 \mathbb P\left(
  \#\{j:s_j(X-zI)\leq W^{-1}\}
  \geq \frac{N^{1+\varepsilon}}W
 \right)=O_{D,\varepsilon,\tau}(N^{-D}).
 \label{eq:historical-count-conjecture}
\end{equation}
\end{conjecture}

This is former Theorem~2.1, now explicitly unproved.  It follows
conditionally from Conjecture~\ref{conj:historical-local-law}:
if \(\lambda_j=s_j(Y_z)^2\) and
\(m_N=N^{-1}\operatorname{Tr}G\), then
\[
 \#\{j:\lambda_j\leq W^{-2}\}
 \leq2N\eta\operatorname{Im}m_N(i\eta,z),
 \qquad \eta=N^{\gamma_0}W^{-2}.
\]
The scalar asymptotic \(|m_c(i\eta,z)|\asymp\eta^{-1/2}\)
and the conjectural local law give an upper bound
\(C N^{1+\gamma_0/2}/W\) with the stated high probability.
Choosing \(\gamma_0\) sufficiently small gives
\eqref{eq:historical-count-conjecture}.  This conditional implication
is not a proof of either conjecture.

\subsection{The incorrect sign and the missing stability estimate}
\label{subsec:historical-proof-error}

The previous argument attempted to adapt resolvent methods from
\cite{MR3068390,MR3230002} to the band setting. We isolate its algebraic
error without assuming that any of the former
probabilistic error estimates have been justified.  Write
\(a=|z|^2\), \(m=m_c(w,z)\),
\(v_i=G_{ii}-m\), and \(\widehat v_i=\mathcal G_{ii}-m\).
The self-consistency equation used in former Proposition~3.13, namely
its equations (3.32)--(3.33), has the form
\begin{equation}
 (m+v_i)^{-1}
 =-w\bigl(1+m+(S^T\widehat v)_i\bigr)
  +\frac a{1+m+(Sv)_i}+e_i.
 \label{eq:historical-correct-start}
\end{equation}
Subtracting \eqref{eq:historical-scalar-equation} and expanding both
reciprocals yields the linear terms
\[
 -\frac{v_i}{m^2}
 =-w(S^T\widehat v)_i-
   \frac a{(1+m)^2}(Sv)_i+\text{remainder}.
\]
Consequently the corrected linearized equation is
\begin{equation}
 v_i=y_1(S^T\widehat v)_i+y_2(Sv)_i+\rho_i,
 \qquad
 y_1=wm^2,\qquad y_2=\frac{am^2}{(1+m)^2}.
 \label{eq:historical-correct-linearization}
\end{equation}
For example, if \(q_i=(Sv)_i\), the exact remainder obtained from
\eqref{eq:historical-correct-start} is
\[
 \rho_i=
 \frac{v_i^2}{m+v_i}
 -\frac{am^2q_i^2}{(1+m)^2(1+m+q_i)}-m^2e_i,
\]
whenever the displayed denominators are nonzero.  No smallness estimate
for this remainder is asserted here.

The previous proof instead used
\(y_2=-am^2/(1+m)^2\).  This is an incorrect sign.  As
\(w=i\eta\to0\), the correct limits are
\begin{equation}
 y_1\longrightarrow-(1-a),\qquad y_2\longrightarrow+a,
 \label{eq:historical-correct-limits}
\end{equation}
whereas \eqref{eq:historical-inverse} assumes that \(y_2\) is close
to \(-a\).  It therefore bounds a different operator and cannot
justify the inversion in the corrected argument.

There is a further reason why changing the sign is not a complete
repair.  The scalar equation gives the exact identity
\begin{equation}
 y_2-y_1=\frac m{1+m}=: \beta,
 \qquad 1-\beta=\frac1{1+m}.
 \label{eq:historical-beta}
\end{equation}
For every doubly stochastic \(S\), the corrected block matrix
\[
 \mathcal B=
 \begin{pmatrix}y_2S&y_1S^T\\y_1S&y_2S^T\end{pmatrix}
\]
has eigenvector \((\mathbf1,-\mathbf1)\) with eigenvalue \(\beta\).
In the symmetric case, subtracting the row and column equations gives
\begin{equation}
 (I-\beta S)(v-\widehat v)=\rho-\widehat\rho.
 \label{eq:historical-difference-channel}
\end{equation}
Here \(|1-\beta|\asymp\sqrt\eta\) for
\(w=i\eta\to0\) and \(0<|z|<1\).

The exact identity
\(\operatorname{Tr}G(w,z)=\operatorname{Tr}\mathcal G(w,z)\)
removes the constant component of \(v-\widehat v\), but it does not
remove all nearly constant modes of a band profile.  For instance, for
the normalized periodic flat band
\[
 s_{ij}=(2W+1)^{-1}\mathbf1\{|i-j|_N\leq W\},\qquad W=o(N),
\]
the first nonconstant Fourier mode \(\phi\) satisfies
\[
 \sum_i\phi_i=0,\qquad S\phi=\lambda\phi,
 \qquad 1-\lambda\asymp(W/N)^2.
\]
On this trace-zero mode, \((I-\beta S)^{-1}\) multiplies by
\((1-\beta\lambda)^{-1}\).  In particular,
\begin{equation}
 |1-\beta\lambda|^{-1}
 \gtrsim\frac1{\sqrt\eta+(W/N)^2}.
 \label{eq:historical-slow-mode}
\end{equation}
This can grow polynomially in \(N\) at the claimed spectral scales
and sublinear polynomial bandwidths.  Thus even after projecting out
the constant mode, the previous polylogarithmic inverse argument is
not justified.

This calculation identifies a failure of the proof, not a
counterexample to the proposed local law.  The actual random residual
in \eqref{eq:historical-difference-channel} can have cancellations
that a worst-case inverse bound does not detect.  A proof at the
former scale would have to control its projections onto the remaining
slow modes, or supply another argument achieving the required
resolvent estimate.  The previous proof supplies neither such an
estimate nor a replacement mechanism.  Double stochasticity selects
the scalar deterministic solution, but it does not by itself give
stability against these nonconstant perturbations.

Accordingly, the former proof of Theorem~3.4 and its use to prove
Theorem~2.1 cannot be retained.  The circular-law conclusions that
depended on them are conjectural at \(W\geq N^{1/2+c}\).
The coarse local law and the \(2/3+c\) result of the present revision
are proved independently, without the historical inverse assumption
\eqref{eq:historical-inverse} and without invoking
Conjectures~\ref{conj:historical-local-law} or
\ref{conj:historical-count}.

% END sections/conjectures
% BEGIN sections/lsv

\section*{Acknowledgment}
The author thanks Giorgio Cipolloni for many insightful discussions on random band matrices.
\section*{Funding}
The author is supported by a Simons Foundation Grant (601948, DJ).
% BEGIN references

% END references
\end{document}